\documentclass[11pt]{llncs}
\usepackage[letterpaper,hmargin=1.1in,vmargin=1.15in]{geometry}
\usepackage{amsmath,amssymb,amsbsy,amsfonts,latexsym,amsopn,amstext,amsxtra,euscript,amscd,color}
\usepackage{array}
\usepackage{booktabs}
\usepackage{graphicx}
\usepackage{slashbox}
\usepackage{psfrag}
\usepackage{epsfig}
\usepackage[sort&compress,numbers,square]{}

\pdfoutput=1

\def\00{{\bf 0}}
\def\+{\oplus}

\def\\{\cr}
\def\({\left(}
\def\){\right)}

\providecommand{\newoperator}[3]{%
  \newcommand*{#1}{\mathop{#2}#3}}

\newoperator{\FD}{\mathrm{FD}}{\nolimits}

\begin{document}
\title{An algorithm based on a new DQM with modified exponential cubic B-splines for solving hyperbolic telegraph
equation in $(2+1)$ dimension}
\author{Brajesh Kumar Singh\thanks{Address for Correspondence: Department of Applied Mathematics, Babasaheb Bhimrao Ambedkar University Lucknow-226 025 (UP) INDIA}, Pramod Kumar}
\institute{Department of Applied Mathematics, School for Physical Sciences, \\ Babasaheb Bhimrao Ambedkar University Lucknow-226 025 (UP)  INDIA \\
\email{bksingh0584@gmail.com, bbaupramod@gmail.com}}
\date{\today}

\maketitle
\begin{abstract}
This paper developed a method called "modified exponential cubic B-Spline differential quadrature (mExp-DQM) for space discretization together with a time integration algorithm" for the numerical computation of hyperbolic telegraph equation in $(2+1)$ dimension. The mExp-DQM is a new differential quadrature method based on modified exponential cubic B-splines as basis which reduces the problem into an amenable system of ordinary differential equations. The resulting system is solved using a time integration algorithm. The stability of the method is also studied by computing the eigenvalues of the coefficients matrices, it is found that the scheme is conditionally stable. The accuracy of the method is illustrated by computing the error between analytical solutions and numerical solutions is measured by using $L_2$ and $L_{\infty}$ error norms for each problem. A comparison of mExp-DQM solutions with the results of the other numerical methods has been carried out for various space sizes and time step sizes.
\end{abstract}


\section{Introduction}
The hyperbolic partial differential equations have a great attention due to its wide range of applications in fields of applied science and engineering, for instance, the hyperbolic partial differential equation models fundamental equations in atomic physics \cite{LS10} and is very useful in understanding various physical phenomena in applied sciences and engineering. It models the vibrations of structures (e.g. buildings, machines and beams). We consider second-order  two-space dimensional linear hyperbolic telegraph equation of the form:
\begin{equation}\label{int-tel-eqn1}
 \frac{\partial^2 u (x, y, t)}{\partial t^2}+2\alpha  \frac{\partial u (x, y, t)}{\partial t} +\beta^2 u(x,y,t)=  \frac{\partial^2 u (x, y, t)}{\partial x^2} + \frac{\partial^2 u (x, y, t)}{\partial y^2}+f(x,y,t),~~ (x,y)\in \Omega, t>0.
\end{equation}
 where $\partial \Omega$ denotes the boundary of the computational domain $\Omega=[0,1]\times [0,1] \subset R^2$ and $\alpha>0, \beta$ are arbitrary constants. Eq. \eqref{int-tel-eqn1} with $\beta=0$ is a damped wave equation while for $\beta >0$ it reduces to telegraph equation. The telegraph equation is more convenient than ordinary diffusion equation in modeling reaction diffusion for such branches of sciences \cite{DG10}, and mostly used in wave propagation of electric signals in a cable transmission line \cite{P90}.

Two space dimensional  initial boundary value problem for second order linear telegraph equation \eqref{int-tel-eqn1} is obtained by combining the equation with the following initial conditions:
\begin{equation}\label{eqn-IC}
 \begin{array}{ll}
 u(x,y,0)=\phi(x,y), \qquad  u_t(x,y,0)=\psi(x,y), \qquad (x, y) \in \Omega,
\end{array}
\end{equation}
and the boundary conditions- Dirichlet boundary condition:
\begin{equation}\label{eqn-DBC}
\left\{ \begin{array}{ll}
 u(0,y,t)=\phi_1(y,t),u(1,y,t)=\phi_2(y,t),\\
 u(x,0,t)=\phi_3(x,t),u(x,1,t)=\phi_4(x,t),
\end{array}  \right. \qquad (x,y)\in \partial \Omega, t>0 \end{equation}
or Neumann boundary conditions:
\begin{equation}\label{eqn-NBC}
\left\{ \begin{array}{ll}
 u_x(0,y,t)=\psi_1(y,t),u(1,y,t)=\psi_2(y,t),\\
 u_y(x,0,t)=\psi_3(x,t),u(x,1,t)=\psi_4(x,t),
\end{array}  \right.  \qquad (x,y)\in \partial \Omega, t>0. \end{equation}
where $\psi, \phi, \psi_i, \phi_i (i=1,2 ,3,4)$ are known smooth functions.

In the recent years, a lot of numerical techniques have been developed for solving hyperbolic telegraph equations in both one and two dimensions, among them: Taylor matrix method \cite{BS11}, dual reciprocity boundary integral method \cite{DG10}, unconditionally stable finite difference scheme \cite{GC07}, implicit difference scheme \cite{M09}, variational iteration method \cite{MYL11},modified B-spline collocation method \cite{MB13}, Chebyshev tau method \cite{SD10}, interpolating scaling function method \cite{LS10},  cubic B-spline collocation method \cite{SR16}.

Two dimensional initial value problem of telegraph equation have been solved by various schemes: Taylor matrix method by B\"{u}lb\"{u}l and Sezer \cite{BS11a} which  converts the telegraph equation into the matrix equation, Two meshless methods- namely meshless local weak-strong
(MLWS) and meshless local Petrov-Galerkin (MLPG) method by Dehghan and Ghesmati \cite{DG10a}, higher order implicit collocation method \cite{DM09}, A polynomial based differential quadrature method \cite{JPM12}, modified cubic B-spline differential quadrature method \cite{MB14}, an unconditionally stable alternating direction implicit scheme \cite{MJ01}, A hybrid method due to Dehghan and Salehi \cite{DS12}, compact finite difference scheme by Ding and Zhang \cite{DZ09} with accuracy of order four in both space and time. Two dimensional linear hyperbolic telegraph equation with variable coefficients has been solved by Dehghan and Shorki \cite{DS09}.

In the recent years, the differential  quadrature method (DQM), developed  by Bellman et al. \cite{BKLV75,BKC72} for the numerical computation of partial differential equations (PDEs), have a great attention among the researchers. After seminal work of Bellman et al.\cite{BKLV75,BKC72}, Quan and Chang \cite{QC89,QC89a}, the DQM has been employed with various types of basis functions, among others, cubic B-spline DQM \cite{KD13,KD13a}, modified cubic B-spline differential quadrature method (MCB-DQM) \cite{AS13,MB14},  DQM based on fourier expansion and Harmonic  function \cite{SC97,SX97,SWB95}, sinc differential quadrature method \cite{KD11a}, generalized DQM \cite{SR92}, polynomial based DQM \cite{K10,JPM12}, quartic B-spline based DQM \cite{BKG15}, Quartic and quintic B-spline methods \cite{KD16}, exponential cubic B-spline DQM \cite{KA15}, extended cubic B-spline DQM \cite{KA16}.

In this paper our aim is to develop modified exponential cubic-B-spline differential quadrature method (mExp-DQM) for hyperbolic partial differential equations. Specially, mExp-DQM is employed for the numerical computation of two dimensional second order linear hyperbolic telegraph equation with both Dirichlet boundary conditions and Neumann boundary conditions. The mExp-DQM is the differential quadrature method based on modified exponential cubic-B-splines as set of new basis functions. The mExp-DQM converts the initial- boundary value system of the telegraph equation  into a initial value system of ODEs, in time. The resulting system of ODEs can be solved by using various time integration algorithm, among them, we prefer SSP-RK54 scheme \cite{GKS09,SR02} due to its reduce storage space, results in less accumulation of the numerical errors. The accuracy and adaptability of the method is illustrated by considering six test problems of the two dimensional telegraphic equations. 

The rest of the paper is organized into five more sections, which follow this introduction. Specifically, Section \ref{sec-disc} deals with
the description of mExp-DQM. Section \ref{sec-imp} is devoted to the procedure for the implementation of mExp-DQM for the problem \eqref{int-tel-eqn1} with the initial conditions \eqref{eqn-IC} and boundary conditions \eqref{eqn-DBC} and \eqref{eqn-NBC}. The stability analysis of the mExp-DQM is studied in Section \ref{sec-stab}. Section \ref{sec-numer} is concerned with the main aim, the numerical study of six test problems, to establish the accuracy of the proposed method in terms of the relative error norm ($R_e$), $L_2$ and $L_\infty$ error norms. Finally, Section \ref{sec-conclu} concludes the paper with reference to critical analysis and research perspectives.

\section{Description of mExp-DQM} \label{sec-disc}
The differential quadrature method is an approximation to derivatives of a function is the weighted sum of the functional values at certain nodes \cite{BKC72}. The weighting coefficients of the derivatives is depend only on grids \cite{SR92}. This is the reason for taking the partitions $P[\Omega]$ of the computational domain $\Omega=\{ (x, y)\in R^2: 0\leq x, y\leq 1\}$ of the problem is distributed uniformly as follows:
 $$P[\Omega] = \{ (x_i, y_j)\in \Omega: h_x=x_{i+1}-x_{i}, h_y=y_{j+1}-y_{j}, i\in \Delta_x, j\in \Delta_y\},$$
where $\Delta_x=\{1,2,\ldots,N_x\}, \Delta_y=\{1,2,\ldots,N_y\}$, and $
h_x=\frac{1}{N_x-1} \mbox{ and } h_y=\frac{1}{N_y-1}$
are the discretization steps in both $x$ and $y$ directions, respectively. That is, a uniform partition in each $x, y$-direction with the following grid points: $$0=x_1 < x_2< \ldots <x_i<\ldots <x_{N_x-1}< x_{N_x}=1,$$
 $$0 = y_1 < y_2< \ldots< y_j< \ldots< y_{N_y-  1}<y_{N_y}=1.$$
Let $(x_i, y_j)$ be the generic grid point and
$$u_{ij}\equiv u_{ij}(t)\equiv u(x_i, y_j, t),~~i\in \Delta_x, j\in \Delta_y.$$

The approximation for $r$-th order derivative of $u(x,y,t)$, for $r \in \{1, 2\}$, with respect to $x, y$ at $(x_i, y_j)$ for $i\in \Delta_{x}, j\in \Delta_{y}$ is given by
\begin{equation}\label{eq-deri1}
\begin{split}
& \frac{\partial^r u}{\partial x^r}(x_i, y_j) = \sum_{\ell=1}^{N_x} a_{i\ell}^{(r)} u_{\ell j}, \qquad  i\in\Delta_x,\\
& \frac{\partial^r u}{\partial y^r}(x_i, y_j) = \sum_{\ell=1}^{N_y} b_{j \ell}^{(r)} u_{i\ell}, \qquad  j\in\Delta_y,
\end{split}
\end{equation}
where the coefficients $a_{ip}^{(r)}$ and $b_{jp}^{(r)}$, the time dependent unknown quantities, are termed as the weighting functions of the $r$th-order derivative, to be calculated using of various type of basis functions.

\noindent The exponential cubic B-splines function $\zeta_i=\zeta_i(x)$ at node $i$ in $x$ direction, reads \cite{KA15,ED15}:
\begin{eqnarray}\label{eq-cbs}
\zeta_i= \frac{1}{h_x^3}\left\{ \begin{array}{ll}
                 b_2\{(x_{i-2}-x)-\frac{1}{p}\sinh(p(x_{i-2}-x))\},  & x \in [x_{i-2}, x_{i-1}) \\
                 a_1+b_1(x_{i}-x)+c_1\exp(p(x_{i}-x))+d_1\exp(p(x_{i}-x)), &   x \in [x_{i-1}, x_{i})\\
                 a_1+b_1(x-x_{i})+c_1\exp(p(x-x_{i}))+d_1\exp(p(x-x_{i})),
                    &  x \in [x_{i}, x_{i + 1})\\
                 b_2\{(x-x_{i+2})-\frac{1}{p}\sinh(p(x-x_{i+2}))\},  &  x \in [x_{i + 1}, x_{i + 2})\\
                 0,  &  \mbox {otherwise}
       \end{array}  \right.
\end{eqnarray}
where
 \begin{equation*}
 \begin{split}
 &a_1=\frac{pc h_x}{pch_x-s}; b_1=\frac{p}{2}\left( \frac{s^2-c(1-c)}{(pch_x-s)(1-c)}\right), b_2=\frac{p}{2(pch_x-s)}, c=\cosh(p h_x), s=\sinh(p h_x), \\
 &c_1=\frac{1}{4} \left\{\frac{\exp(-ph_x)(1-c)+s(\exp(-ph_x)-1) }{(pch_x-s)(1-c)}\right\}, d_1=\frac{1}{4} \left\{\frac{\exp(ph_x)(c-1)+s(\exp(ph_x)-1) }{(pch_x-s)(1-c)}
 \right\}.
 \end{split}
\end{equation*}
 The set $\{\zeta_{0}, \zeta_{1}, \zeta_{2}, \ldots, \zeta_{N_x}, \zeta_{N_x+1}\}$ forms a basis over the interval $[a, b]$. The values of $\zeta_i$ and its first and second derivatives in the grid point $x_j$, denoted by $\zeta_{ij}:=\zeta_i(x_j)$, $\zeta'_{ij}:=\zeta'_i(x_j)$ and $\zeta''_{ij}:=\zeta''_i(x_j)$, respectively,  read:
\begin{eqnarray}\label{tab-coeff}
       \zeta_{ij}= \left\{ \begin{array}{ll}
                 1,   &\mbox{ if }  i-j=0 \\
                 \frac{s-ph}{2(pc h_x-s)},  & \mbox{ if } i-j=\pm 1\\
                 0,  & \mbox {otherwise}
       \end{array}\right.\\
       \zeta'_{ij}= \left\{ \begin{array}{ll}
                  - \frac{p(1-c)}{2(pc h_x - s)},  & \mbox{ if }  i-j=1\\
                   \frac{p(1-c)}{2(pc h_x - s)},  & \mbox{ if }  i-j=- 1\\
                 0, &\mbox {otherwise}
       \end{array}\right.\mbox{~}\\
       \zeta''_{ij}= \left\{ \begin{array}{ll}
                -\frac{p^2 s}{(pc h_x - s)},  &\mbox{ if }  i-j=0 \\
                 \frac{p^2 s}{2(pc h_x - s)},  &\mbox{ if }  i-j=\pm 1\\
                 0  &  \mbox {otherwise}
       \end{array}\right.
\end{eqnarray}
The modified exponential cubic B-splines basis functions are obtained by modifying the exponential cubic B-spline basis function \eqref{eq-cbs} as follows \cite{AS13}:
\begin{flushleft}
\begin{equation}\label{eq-modi-CBs}
\left\{ \begin{split}
& \psi_1 (x) = \zeta_1(x) + 2\zeta_0(x)\\
&  \psi_2 (x) = \zeta_2(x)  - \zeta_0(x)\\
& \vdots\\
& \psi_j (x) = \zeta_j(x), \mbox{ for } j = 3,4, \ldots, N_x-2\\
& \vdots\\
&  \psi_{N_x-1}(x) = \zeta_{N_x-1}(x)  - \zeta_{N_x+1}(x)\\
&  \psi_{N_x}(x) = \zeta_{N_x}(x)  + 2 \zeta_{N_x+1}(x)
\end{split}\right.
\end{equation}
\end{flushleft}
The set $\{\psi_1, \psi_2,\ldots, \psi_{N_x}\}$ is a basis over the set $[a, b]$. Analogously procedure is followed for $y$ direction.

\subsection{\textbf{The evaluation of the weighting coefficients $a_{ij}^{(r)}$ and  $ b_{ij}^{(r)} (r=1,2)$}}
In order to evaluate the weighting coefficients $a_{ip}^{(1)}$ of first order partial derivative in Eq. \eqref{eq-deri1}, we use the modified exponential cubic B-spline $\psi_p(x)$, $p\in\Delta_x$ in DQ method as set of basis functions. Write $\psi'_{pi}:=\psi_p'(x_i)$ and $\psi_{p\ell}:=\psi_p(x_\ell)$, In DQ method, the approximate values of the first-order derivative is obtained as:
\begin{equation}\label{eq-deri1approxx}
\psi'_{pi} = \sum_{\ell=1}^{N_x} a_{i\ell}^{(1)}  \psi_{p\ell}, \qquad  p,i \in \Delta_x.
\end{equation}
Setting  $\Psi=[\psi_{p\ell}]$, $A=[a_{i\ell}^{(1)}]$ (the unknown weighting coefficient matrix), and $\Psi'=[\psi'_{p\ell}]$, then Eq. \eqref{eq-deri1approxx} can be re-written as the following set of system of linear equations:
\begin{equation}\label{eq-tri_syst}
\Psi A^T=\Psi'.
\end{equation}
Let $\omega = \frac{p(1-c)h_x}{pc h_x-s}$ and $\theta =\frac{s-ph_x}{2(pc h_x-s)}$, then the coefficient matrix $\Psi$ of order $N_x$ can be obtained from  \eqref{tab-coeff} and \eqref{eq-modi-CBs}:
\begin{equation*}
 \Psi= \left[
  \begin{array}{cccccccc}
  \omega  $\quad$ & \theta $\quad$&    $\quad$&           $\quad$&           $\quad$&         $\quad$&         \\
  $0    \quad$& 1 $\quad$& \theta          $\quad$&           $\quad$&           $\quad$&         $\quad$&         \\
       $\quad$& \theta   $\quad$& 1       $\quad$& \theta      $\quad$&           $\quad$&         $\quad$&         \\
       $\quad$&       $\quad$&  \ddots      $\quad$&   \ddots  $\quad$&   \ddots  $\quad$&         $\quad$&         \\
       $\quad$&       $\quad$&              $\quad$&   \theta     $\quad$&   1 $\quad$&    \theta  $\quad$&          \\
       $\quad$&       $\quad$&              $\quad$&           $\quad$&   \theta     $\quad$&   1$\quad$&    $0$   \\
       $\quad$&       $\quad$&              $\quad$&           $\quad$&           $\quad$&   \theta   $\quad$&   \omega   \\
       \end{array}
\right]
\end{equation*}
and in particular the columns of the matrix $\Psi'$ read:
\begin{equation*}
\Psi'[1] = \left[ \begin{array}{c}
    \omega/h_x\\
    -\omega/h_x\\
    0 \\
  \vdots  \\
          \\
  0     \\
  0     \\
\end{array}
\right],
\Psi'[2] =\left[\begin{array}{c}
    \omega/2h_x \\
    $0$    \\
    -\omega/2h_x  \\
     $0$   \\
   \vdots  \\
           \\
    $0$    \\
\end{array}\right],
\ldots,
\Psi'[N_x-1] =  \left[ \begin{array}{c}
     $0$     \\
     \vdots  \\
             \\
     $0$     \\
    \omega/2h_x  \\
     $0$     \\
     -\omega/2h_x   \\
\end{array}
\right], \mbox{ and }
\Psi'[N_x] =\left[ \begin{array}{c}
    $0$ \\
        \\
\vdots  \\
        \\
$0$     \\
 \omega/h_x \\
-\omega/h_x   \\
\end{array}
\right].
\end{equation*}
It is remarked that the exponential cubic B-splines are modified in order to have a diagonally dominant coefficient matrix $\Psi$, see Eq. \eqref{eq-tri_syst}. To calculate the weighting coefficients, we solve the system \eqref{eq-tri_syst} using the well known Thomas Algorithm \cite{Lee}. Similarly, the weighting coefficients $b_{i \ell}^{(1)}$ can be evaluated by considering the grids in the $y$ direction.

In similar manner, the weighting coefficients $a_{ip}^{(r)}$ and $b_{ip}^{(r)}$, for $r\ge 2$, can be calculated using the weighting functions in quadrature formula for second order derivative on the given basis functions. But, in the present paper, we prefer the following recursive formulae \cite{SR92}:
\begin{flushleft}
 \begin{equation}\label{eq-coeff2}
\left\{ \begin{split}
& a_{i j}^{(r)} = r \left( a_{i j}^{(1)} a_{i i}^{(r-1)} - \frac{a_{i j}^{(r-1)}}{ x_i - x_j}\right), i\ne j: i, j\in \Delta_x,
\\&  a_{i i}^{(r)} = - \sum_{i = 1, i \ne j}^{N_x} a_{i j}^{(r)}, i= j: i, j\in \Delta_x.\\
& b_{i j}^{(r)} = r \left(b_{i j}^{(1)} b_{i i}^{(r-1)} - \frac{b_{i j}^{(r-1)}}{ y_i - y_j}\right), i\ne j: i, j \in \Delta_y
 \\& b_{i i}^{(r)} = - \sum_{i = 1, i \ne j}^{N_y} b_{i j}^{(r)}, i = j: i, j \in \Delta_y.\\
\end{split}\right.
\end{equation}
\end{flushleft}
\section{The mExp-DQM for the telegraph equation} \label{sec-imp}
First, we set $u_t=v$ and thus $u_{tt}=v_t$. Keeping all above in mind, the second order telegraph equation \eqref{int-tel-eqn1} with the initial condition transforms to initial valued coupled system of first order differential equations as follows:
\begin{equation}\label{trans-tel-eqn2}
\left\{ \begin{array}{ll}
\frac{\partial u (x, y, t)}{\partial t}=v(x,y,t)\\ \\
\frac{\partial v (x, y, t)}{\partial t}= - 2\alpha  v (x, y, t) - \beta^2 u(x,y,t)+ \frac{\partial^2 u (x, y, t)}{\partial x^2} + \frac{\partial^2 u (x, y, t)}{\partial y^2}+f(x,y,t), (x,y)\in \Omega, t>0, \\
 u(x,y,0)=\phi(x,y), \qquad  v(x,y,0)=\psi(x,y), \qquad (x, y) \in \Omega.
\end{array}  \right.\end{equation}
Further, on setting $f(x_i, y_j,t)=f_{ij}$, the mExp-DQM transforms Equation \eqref{trans-tel-eqn2} to
\begin{equation}\label{trans-tel-eqn3}
\left\{ \begin{split}
& \frac{\partial u_{ij}}{\partial t}=v_{ij}\\
& \frac{\partial v_{ij}}{\partial t}= \sum_{\ell=1}^{N_x}a_{i\ell}^{(2)}u_{\ell j}+\sum_{\ell=1}^{N_y}b_{j\ell}^{(2)}u_{i\ell}-2\alpha v_{ij}-\beta^2 u_{ij}+f_{ij},\\
&  u_{ij}(t=0)=\phi_{ij}, \quad
 v_{ij}(t=0)=\psi_{ij}, \quad i \in \Delta_x, ~~  j\in \Delta_y.
\end{split}\right.
\end{equation}

Next, further simplification is not required in case of Dirichlet boundary conditions. In this case the solution on the boundary can be read directly from the conditions \eqref{eqn-DBC} as:
\begin{flushleft}
 \begin{equation}\label{eq-DBC-dis}
\left. \begin{split}
& u_{1j} = \phi_{1}(y_j, t)=\phi_{1}(j), \qquad u_{N_xj} = \phi_{2}(y_j, t)=\phi_{2}(j), \qquad  j\in \Delta_y, \\
& u_{i1} = \phi_{3}(x_i, t)=\phi_{3}(i), \qquad u_{iN_y} = \phi_{4}(x_i, t)=\phi_{4}(i), \qquad i \in \Delta_x,
\end{split}  \right\} t\geq0.
\end{equation}
\end{flushleft}
But, for Neumann or mixed type boundary conditions, further simplification is required. In this case, the solutions at the boundary are obtained by using mExp-DQM on the boundary. This yields a system of linear equations, and on solving resulting system of the linear, the desired solution is obtained at the boundaries.

In particular, Eq. \eqref{eq-deri1} with $r=1$ and the Neumann boundary conditions \eqref{eqn-NBC} at $x= 0$ and $x=1$ can be written as
\begin{flushleft}
 \begin{equation}\label{eq-N1-dis}
 \begin{split}
& \sum_{\ell=1}^{N_x} a_{1\ell}^{(1)} u_{\ell j} = \psi_1(j), \qquad
 \sum_{\ell=1}^{N_x} a_{N_x \ell }^{(1)} u_{\ell j} = \psi_2(j),
\end{split} \qquad j\in \Delta_y.
\end{equation}
\end{flushleft}
In terms of matrix system for $u_{1j}, u_{N_x j}$, the above equation can be rewritten as
\begin{equation}\label{eq-N2-dis}
\left[\begin{array}{cc}
a_{11}^{(1)}  & a_{1N_x}^{(1)}  \\
a_{N_x1}^{(1)} & a_{N_x N_x}^{(1)} \\
\end{array} \right] \left[
\begin{array}{c}
u_{1j} \\
u_{N_xj} \\
\end{array} \right]
=
\left[
\begin{array}{c}
S_1(j)\\
S_2(j) \\
\end{array} \right],
\end{equation}
where $S_1(j)= \psi_1(j)-\sum_{\ell =2}^{N_x-1} a_{1\ell}^{(1)} u_{\ell j}$ and $ S_2(j)= \psi_1(j)- \sum_{\ell=2}^{N_x-1} a_{N_x \ell}^{(1)} u_{\ell j}$.

\noindent On solving \eqref{eq-N2-dis}, we have
\begin{equation}\label{eq-N3-dis}
u_{1j} = \frac{S_1(j) a_{N_xN_x}^{(1)} - S_2(j) a_{1N_x}^{(1)} }{a_{11}^{(1)}a_{N_xN_x}^{(1)}-a_{N_x 1}^{(1)}a_{1N_x}^{(1)} },  \qquad
u_{N_x j} = \frac{S_2(j) a_{11}^{(1)} - S_1(j) a_{N_x 1}^{(1)}}{a_{11}^{(1)}a_{N_xN_x}^{(1)}-a_{N_x 1}^{(1)}a_{1N_x}^{(1)} }, \qquad j\in \Delta_y.
\end{equation}

Similarly, Eq. \eqref{eq-deri1} with $r=1$ and Neumann boundary conditions \eqref{eqn-NBC} at $y= 0$ and $y=1$ can be written as
\begin{equation*}
 b_{11}^{(1)} u_{i1} + b_{1N_y}^{(1)} u_{iN_y} = S_3(i),  \quad 
 b_{N_x 1}^{(1)} u_{i1} + b_{N_y N_y}^{(1)} u_{iN_y}     =  S_4(i), \qquad i\in \Delta_x.
\end{equation*}
where $S_3(i)= \psi_3(i) - \sum_{\ell=2}^{N_y-1} b_{1\ell }^{(1)} u_{i\ell }$ and $ S_4(i)= \psi_4(i) - \sum_{\ell=2}^{N_y-1} b_{N_y \ell}^{(1)} u_{i\ell}$.

\noindent On solving the above system of equations for the boundary values $u_{i1} $ and $ u_{iN_y}$, we have
\begin{equation}\label{eq-N2y-dis}
 \begin{split}
&u_{i1} = \frac{S_3(i) b_{N_yN_y}^{(1)} - S_4(i) b_{1N_y}^{(1)} }{b_{11}^{(1)}b_{N_y N_y}^{(1)}-b_{N_y 1}^{(1)}b_{1N_y}^{(1)} }, \qquad  u_{iN_y} = \frac{S_4(i) b_{11}^{(1)} - S_3(i) b_{N_y 1}^{(1)}}{b_{11}^{(1)}b_{N_y N_y}^{(1)}-b_{N_y 1}^{(1)}b_{1 N_y}^{(1)} },
\end{split}
\qquad i\in \Delta_x.
\end{equation}

Finally, on using boundary values $u_{1j}, u_{N_x j}, u_{i1}$ and $u_{iN_y}$ obtained from either boundary conditions (Dirichlet boundary conditions \eqref{eqn-DBC} or Neumann boundary conditions  \eqref{eqn-NBC}), Eq. \eqref{trans-tel-eqn3} can be rewritten as follows:
\begin{equation}\label{eq-ode-finl1}
 \left\{ \begin{split}
& \frac{\partial u_{ij}}{\partial t}=v_{ij}\\
& \frac{\partial v_{ij}}{\partial t}= \sum_{\ell=2}^{N_x-1}a_{i\ell}^{(2)}u_{\ell j}+\sum_{\ell=2}^{N_y-1}b_{j\ell}^{(2)}u_{i\ell}-2\alpha v_{ij}-\beta^2 u_{ij}+K_{ij},\\
&  u_{ij}(t=0)=\phi_{ij}, \quad
 v_{ij}(t=0)=\psi_{ij},
\end{split}\right.
\end{equation}
where $2\leq i\leq N_x-1, 2\leq j\leq N_y-1$ and
 \begin{equation}\label{eq-F}
 \begin{split}
&K_{ij}=f_{ij}+ a_{i1}^{(2)} u_{1j} + a_{i N_x}^{(2)} u_{N_x j} +b_{j1}^{(2)} u_{i1} + b_{j N_y}^{(2)}u_{iN_y} .
\end{split}
\end{equation}

 A lot time integration schemes have been proposed for the numerical computation of initial valued system of differential equations \eqref{eq-ode-finl1}, among others, the SSP-RK scheme allows low storage and large region of absolute property \cite{GKS09,SR02}. In what follows, we adopt SSP-RK54 scheme, is strongly stable for nonlinear hyperbolic differential equations:
\begin{equation*}
\begin{split}
&u^{(1)} = u^m+ 0.391752226571890 \triangle t L(u^m) \\
&u^{(2)} = 0.444370493651235  v^m+ 0.555629506348765 u^{(1)} + 0.368410593050371 \triangle t L(u^{(1)}) \\
&u^{(3)}= 0.620101851488403 u^m+ 0.379898148511597 u^{(2)} + 0.251891774271694 \triangle t L(u^{(2)}) \\
&u^{(4)}= 0.178079954393132 u^m+ 0.821920045606868 u^{(3)} + 0.544974750228521 \triangle t L(u^{(3)}) \\
&u^{m + 1}= 0.517231671970585 u^{(2)}+ 0.096059710526147 u^{(3)} \\
&~~ + 0.063692468666290 \triangle t L(u^{(3)}) + 0.386708617503269 u^{(4)}+ 0.226007483236906 \triangle t L(u^{(4)})
\end{split}
\end{equation*}

\section{Stability analysis} \label{sec-stab}
In compact form , the system \eqref{eq-ode-finl1} can be rewritten as follows:
 \begin{equation}\label{eq-TEL-ode1}
 \left\{ \begin{split}
&\frac{dU}{dt} = AU  + G, \\
& U (t=0)=U_0
\end{split}\right.
\end{equation}
where
\begin{enumerate}
   \item [$1)$] $A =\left[\begin{array}{ccc}
     O \qquad &  &\qquad I \\
     B \qquad & &\qquad-2\alpha I\\
       \end{array} \right],$ \qquad $G = \left[\begin{array}{c}
     O_1 \\
     K \\
       \end{array} \right]$, $U = \left[\begin{array}{c}
     u \\
     v\\
       \end{array} \right], $ and $U_0 = \left[\begin{array}{c}
     \phi \\
     \psi\\
       \end{array} \right]$
   \item [$1)$] $O$ and $O_1$ are null matrices;
   \item [$2)$] $I$ is the identity matrix of order $(N_x-2)(N_y-2)$;
   \item [$3)$] $U=(u,v)^T$ the vector solution at the grid points:

 \noindent  $u=(u_{22},u_{23}, \ldots,u_{2(N_y-1)}, u_{32},u_{33}, \ldots,u_{3(N_y-1)}, \ldots, u_{(N_x-1)2}, \ldots,u_{(N_x-1)(N_y-1)})$.

 \noindent  $v=(v_{22},v_{23}, \ldots,v_{2(N_y-1)}, v_{32},v_{33}, \ldots,v_{3(N_y-1)},  \ldots, v_{(N_x-1)2}, \ldots,v_{(N_x-1)(N_y-1)})$.

   \item [$4)$] $K=(K_{22},K_{23}, \ldots,K_{2(N_y-1)}, K_{32}, \ldots,K_{3(N_y-1)}, \ldots K_{(N_x-1) 2}, \ldots K_{(N_x-1) (N_y-1)} $, where $K_{ij}$, for $i \in \Delta_x, j \in \Delta_y$  is calculated from Eq. \eqref{eq-F}.
  \item [$5)$] $B = -\beta^2 I + B_x + B_y$, where $B_{x}$ and $B_{y}$ are the following matrices (of order $(N_x-2)(N_y-2)$) of the weighting coefficients $a_{ij}^{(2)}$ and $b_{ij}^{(2)}$:
\begin{equation}
\begin{array}{ll}
 B_x = \left[
    \begin{array}{cccc}
      a_{22}^{(2)} I_x & a_{23}^{(2)} I_x & \ldots & a_{2(N_x-1)}^{(2)} I_x \\
      a_{32}^{(2)} I_x & a_{33}^{(2)} I_x &\ldots  & a_{3(N_x-1)}^{(2)} I_x \\
      \vdots & \vdots & \ddots & \vdots \\
      a_{(N_x-1)2}^{(2)} I_x & a_{(N_x-2)3}^{(2)} I_x & \ldots  & a_{(N_x-1)(N_x-1)}^{(2)} I_x \\
    \end{array}
  \right], ~~ & \begin{array}{ll}
   B_y = \left[
    \begin{array}{cccc}
      M_y & O_y & \ldots & O_y \\
      O_y & M_y &\ldots  & O_y \\
      \vdots & \vdots & \ddots & \vdots \\
      O_y & O_y & \ldots  & M_y \\
    \end{array}
  \right]
\end{array}
\end{array}
\end{equation}
where  identity matrix, $I_x$, and null matrix, $O_y$, both are of order $(N_y-2)$ and
\begin{equation*}
\begin{array}{ll}
M_y = \left[
    \begin{array}{cccc}
      b_{22}^{(2)} & b_{23}^{(2)} & \ldots & b_{2(N_y-1)}^{(2)} \\
      b_{32}^{(2)} & b_{33}^{(2)} & \ldots & b_{3(N_y-1)}^{(2)} \\
      \vdots & \vdots & \ddots & \vdots \\
      b_{(N_y-1)2}^{(2)} & b_{(N_y-1)3}^{(2)} & \ldots & b_{(N_y-1)(N_y-1)}^{(2)}
    \end{array}
  \right]
\end{array}
\end{equation*}
 \end{enumerate}

\begin{figure}
\centering
\includegraphics[height=6.0cm, width=10.95cm]{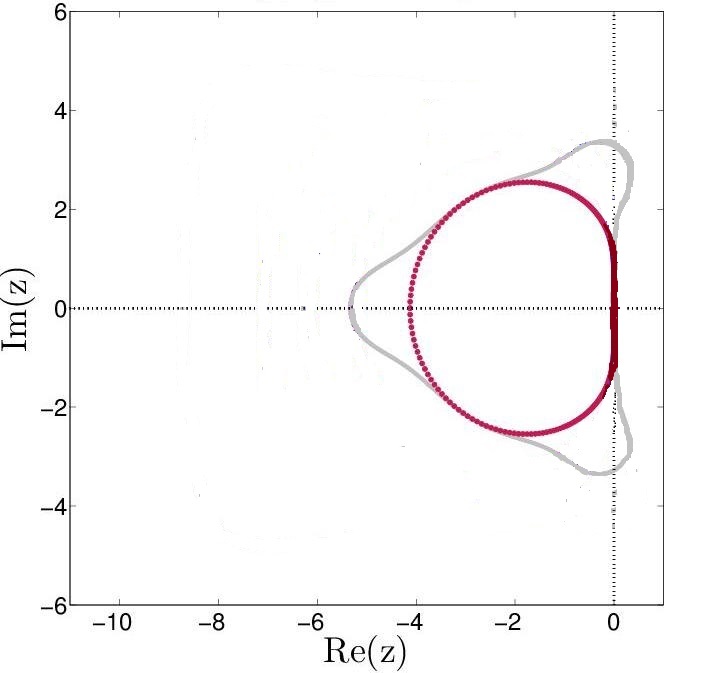}
\caption{ Stability region for SSP-RK54 scheme}\label{Stab-R}
\includegraphics[height=7.0cm, width=14.95cm]{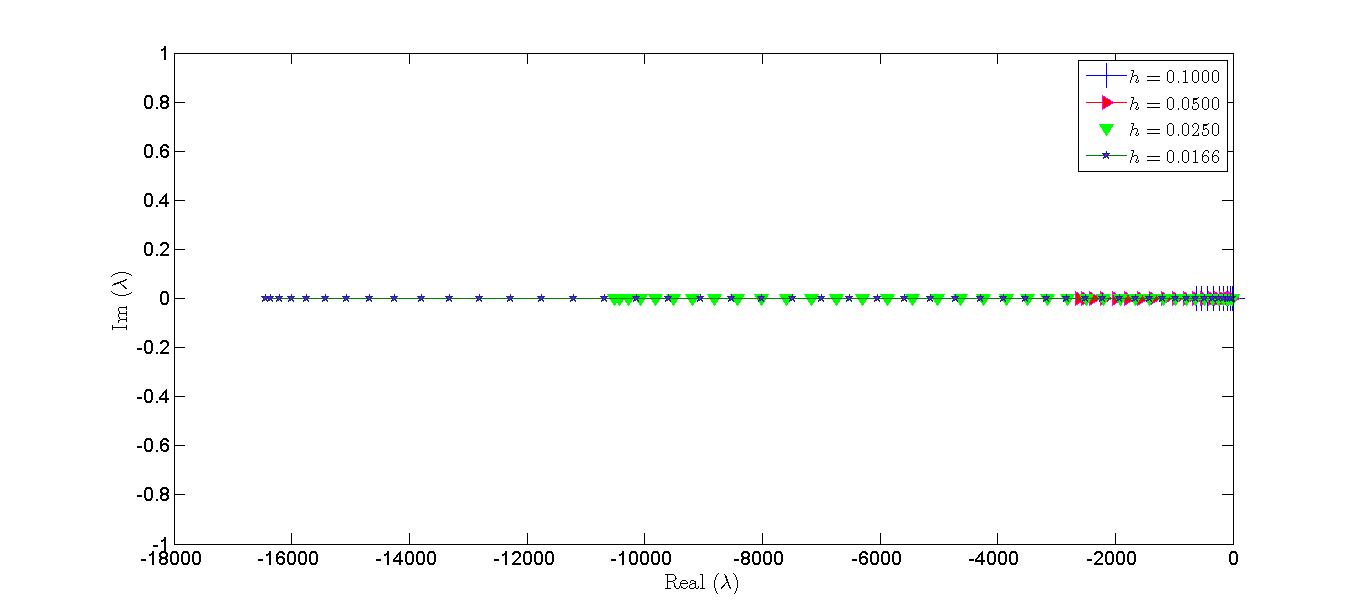}
\caption{ Eigenvalues of $B_x$ and $B_y$ for different grid sizes $h=0.1, 0.01, 0.025, 0.016$}\label{EigB12}
\end{figure}

\noindent The stability of mExp-DQM for the telegraph equation \eqref{int-tel-eqn1} depends on the stability of the system of ODEs defined in \eqref{eq-TEL-ode1}. It is to be noticed that whenever the system of ODEs \eqref{eq-TEL-ode1} is unstable, the proposed method for temporal discretization may not converge to the exact solution. Moreover, being the exact solution can directly obtained by means of the eigenvalues method, the stability of \eqref{eq-TEL-ode1} depends on the eigenvalues of the coefficient matrix $A$ \cite{Jain83}.  In fact, the stability region is the set
$\mathcal{S}= \{z \in C: \mid R(z)\mid \leq 1, z = \lambda_A \triangle t \}$, where $R(.)$ is the stability function and $\lambda_A$ is the eigenvalue of the coefficient matrix $A$. For SSP-RK54 scheme the stability region is depicted in Fig \ref{Stab-R}, see \cite[Fig. 5]{KBYD}. This evident that the sufficient condition for the stability of the system \eqref{eq-TEL-ode1} is that to each eigenvalue $\lambda_A$ of the coefficient matrix A, $\lambda_A \triangle t \in \mathcal{S}$, and hence, the real part of each eigenvalue is necessarily either zero or negative.

Let $\lambda_{A}$ be an eigenvalue of $A$ associated with the eigenvector $(X_1,X_2)^T$, where each component is a vector of order $(N_x-2)(N_y-2)$. Then from Eq. \eqref{eq-TEL-ode1} we have
\begin{equation}\label{eq-egen}
A \left[\begin{array}{c}
     X_1 \\
     X_2\\
       \end{array} \right]=\left[\begin{array}{ccc}
     O \qquad & &\qquad  I \\
     B \qquad & &  \qquad  -2 \alpha I\\
       \end{array} \right]
\left[\begin{array}{c}
     X_1 \\
     X_2\\
       \end{array} \right] = \lambda_A \left[\begin{array}{c}
     X_1 \\
     X_2\\
       \end{array} \right],
\end{equation}
 which yields
\begin{equation}\label{eqn-eign-a}
IX_2 = \lambda_A X_1,
\end{equation}
and
\begin{equation}\label{eqn-eign-b}
BX_1 - 2\alpha X_2 = \lambda_A X_2.
\end{equation}
Simplifying Eq. \eqref{eqn-eign-a} and Eq. \eqref{eqn-eign-b}, we get
\begin{equation}\label{eq-eigenB}
 BX_1 = \lambda_A (\lambda_A + 2\alpha) X_1.
 \end{equation}
which confirms that the eigenvalue $\lambda_B$ of $B$ is $\lambda_B=\lambda_A (\lambda_A + 2 \alpha)$. By definition the matrix $B$ is:
\begin{equation}\label{eq-matrixB}
B= -\beta^2 I + B_x + B_y,
 \end{equation}

Now, we compute the eigenvalues $\lambda$ of $B_x + B_y$ for $p=1$ and  different grid sizes: $h_x=h_y=h=0.1, 0.01, 0.025, 0.016$, and depicted in Fig \ref{EigB12}.

It is evident from Eq. \eqref{eq-matrixB} and Fig. \ref{EigB12} that for different values of the grid sizes the computed eigenvalues $\lambda_B=\lambda-\beta^2$ of $B$ are real negative numbers, that is
\begin{equation}\label{eq-A}
Re\left(\lambda_B\right) \leq 0 \mbox{~~  and ~~} Im\left(\lambda_B\right) = 0,
\end{equation}
where $Re(z)$ and $Im(z)$ denote the real and the imaginary part of $z$, respectively.

Let $\lambda_A = x+\iota y$, then
\begin{equation}\label{eq-B}
\begin{split}
\lambda_B &=\lambda_A (\lambda_A + 2\alpha)\\
& =x^2-y^2 + 2 \alpha x + 2 \iota (x+\alpha)y.
\end{split}
\end{equation}
According to Eq. \eqref{eq-A} and Eq. \eqref{eq-B}, we have
\begin{equation}\label{eq-sta1}
 \left\{ \begin{split}
& x^2-y^2 + 2 \alpha x < 0\\
&(x+\alpha)y=0\\
\end{split}\right.
\end{equation}
The possible solutions of Eq.\eqref{eq-sta1} are
\begin{enumerate}
   \item  [$1)$] If $y\ne 0$, then $x=-\alpha,$
   \item  [$2)$] If $y=0$, then $(x+\alpha)^2 < \alpha^2.$
 \end{enumerate}
 In each case $x$ is negative whenever $\alpha>0$. One a given grid size, one can find a sufficient small value of $\Delta t$ so that $\Delta t \lambda_A$, for each eigenvalue $\lambda_A$ of matrix $A$,  lie inside the stability region $\mathcal{S}$ of SSP-RK54 scheme.  This shows that the mExp-DQM produces stable solutions for two dimensional second order linear telegraph equation.

\section{Numerical experiments and discussion} \label{sec-numer}
This section with the main goal of the paper, the computation of numerical solution of the telegraph equation. The accuracy and the efficiency of this method is measured for six numerical examples in terms of the discrete relative error $R_e$, $L_2$ and $L_{\infty}$ error norms:
\begin{equation*}
 L_2 = \(h \sum_{j=1}^{N} \left[u_j^{exact} - u_j^*\right]^2 \)^{1/2}, \mbox{ and } L_\infty = \max_{j=1}^{N} \left|u_j^{exact} - u_j^*\right|,
\end{equation*}
where $u_j^*$ represent the numerical solution at node $j$. Throughout the section, we have taken equal grid size in each direction, i.e., $h_x=h_y=h$.

\begin{figure}
\centering
\includegraphics[height=7.0cm,width=7.6cm]{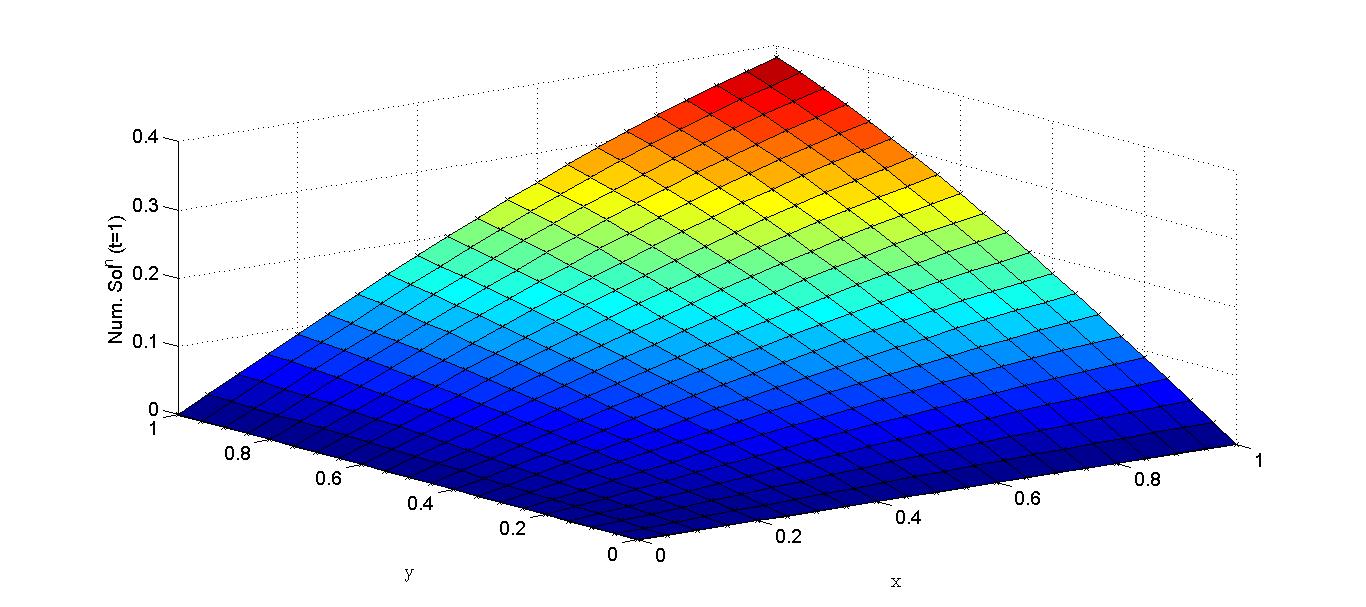}
\includegraphics[height=7.0cm,width=7.6cm]{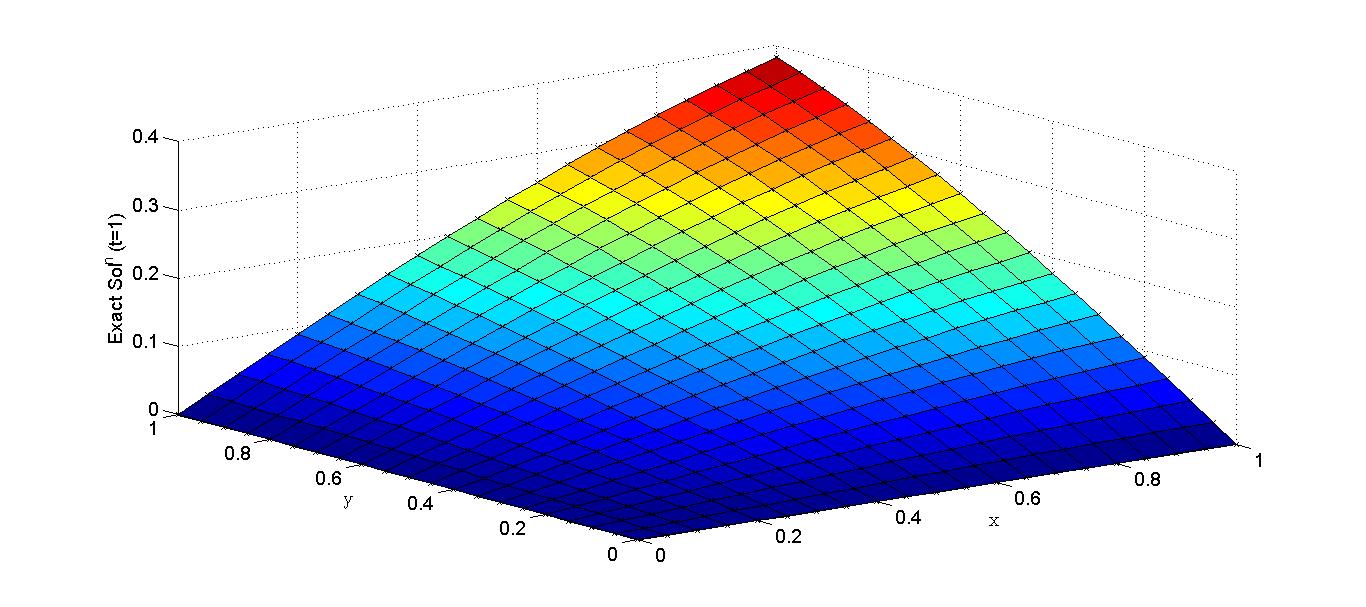}
\includegraphics[height=7.0cm,width=7.6cm]{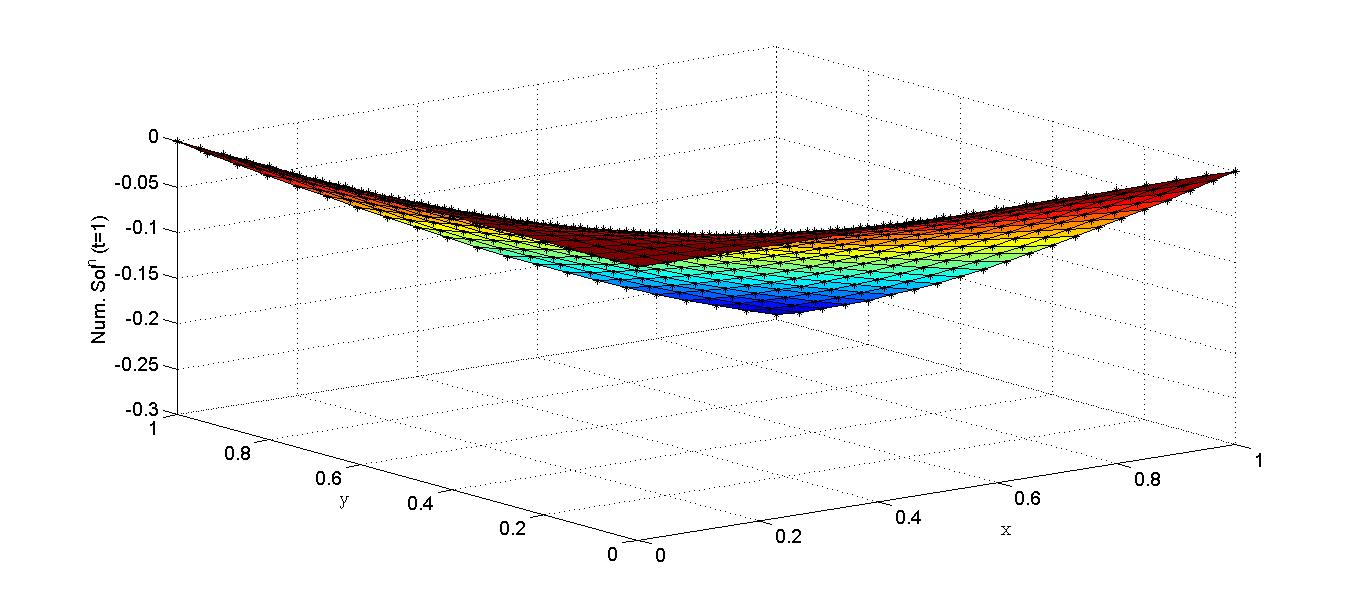}
\includegraphics[height=7.0cm,width=7.6cm]{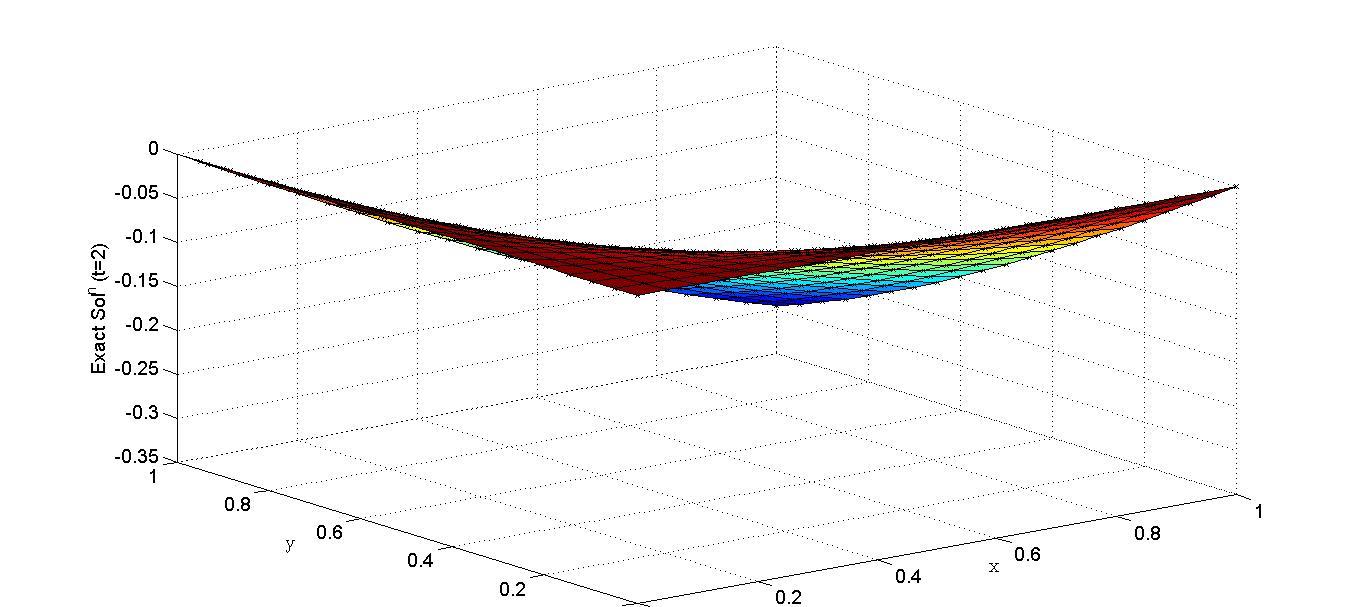}
\includegraphics[height=7.0cm,width=7.6cm]{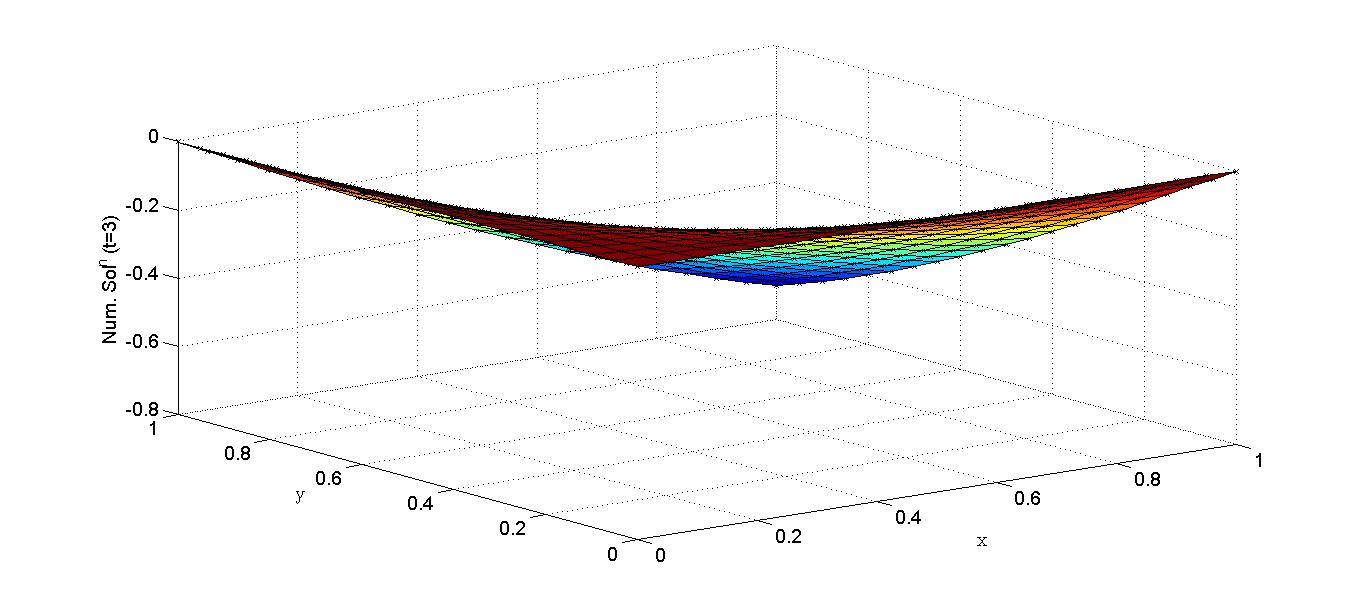}
\includegraphics[height=7.0cm,width=7.6cm]{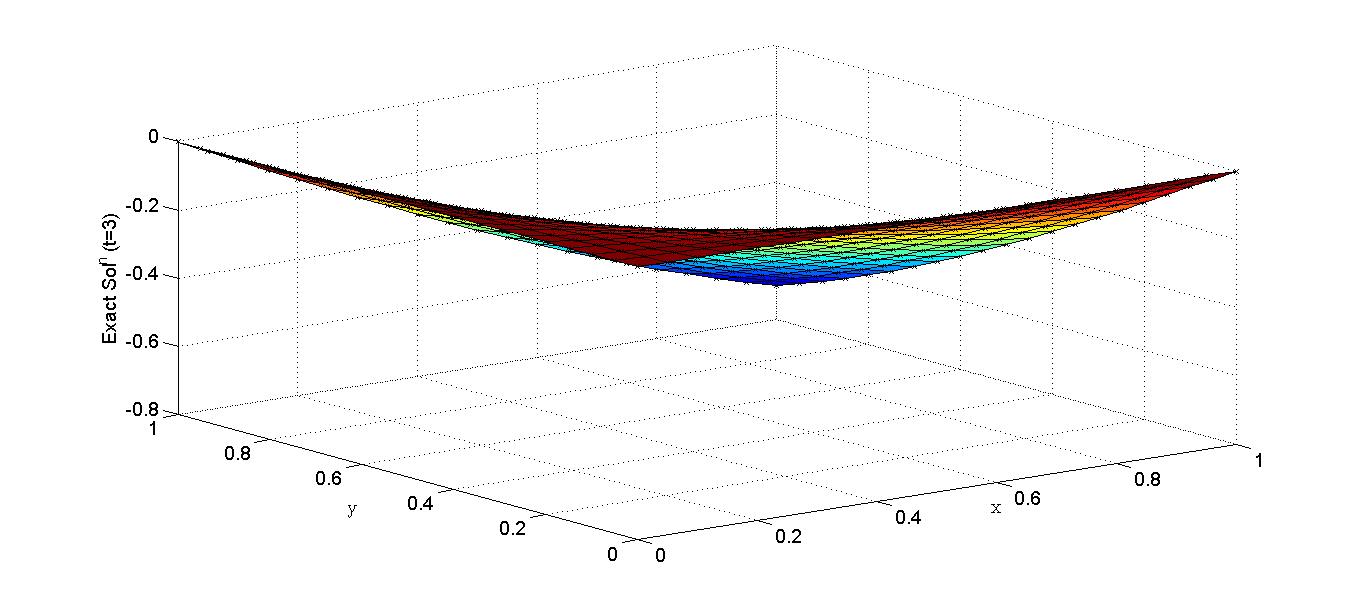}
\caption{Plots of mExp-DQM solution and exact solution of Example \ref{ex1}}\label{fig1.1}
\end{figure}

\begin{example}\label{ex1}
Consider the telegraph equation \eqref{int-tel-eqn1} in the region $\Omega$ with $\alpha=\beta=1$, $f(x,y,t)=2(\cos t-\sin t)\sin x \sin y$,
$ \phi(x,y)=\sin x. \sin y; \psi(x,y)=0$, and the Dirichlet boundary conditions:
\begin{equation}\label{DC-ex1}
\left\{ \begin{array}{ll}
 \phi_1(y,t)=0, \phi_2(y,t)=\cos t \sin(1)\sin y, & 0\le y\le 1, \\ \\
  \phi_3(x,t)=0,   \phi_4(x,t)=\cos t \sin x \sin (1), \qquad& 0\le x \le 1,
\end{array}  \right.\end{equation}
The exact solution \cite{DG10} is
\begin{equation}\label{ex1-exact-soln}
 u(x,y,t)=\cos t \sin x \sin y
\end{equation}
The computed relative error ($R_e$), $L_2, L_\infty$ error norms are compared with the recent results of Mittal and Bhatia \cite{MB14} at different time levels $t\leq 10$, reported in Table \ref{tab-ex1.1} and Table \ref{tab-ex1.2} with the parameters $\Delta t=0.01, h=.1, p=1$ and $\Delta t=0.001$ and $h=.05, p=0.15, 1$, respectively.

The comparison of computed physical solution behavior with the the exact solution behavior at $t=1,2,3$ is depicted in Fig. \ref{fig1.1} with $\Delta t=0.001, h=.05, p=1$. The findings shows that the proposed solution are much better than that of Mittal and Bhatia\cite{MB14}, and are in excellent agreement with the exact solutions. The computation time is slightly more than Mittal and Bhatia\cite{MB14} for large $t$ due to selection of SSP-RK54 scheme instead of SSP-RK43 scheme in time integration.

\end{example}

\begin{figure}
\centering
\includegraphics[height=7.0cm,width=7.6cm]{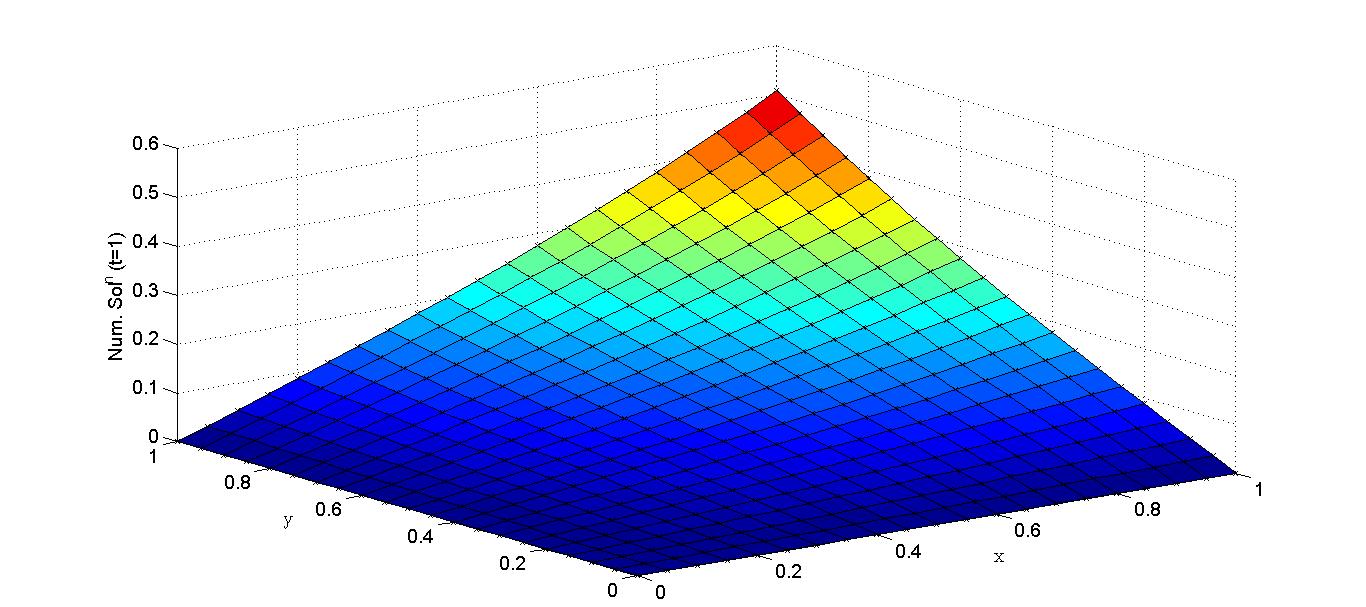}
\includegraphics[height=7.0cm,width=7.6cm]{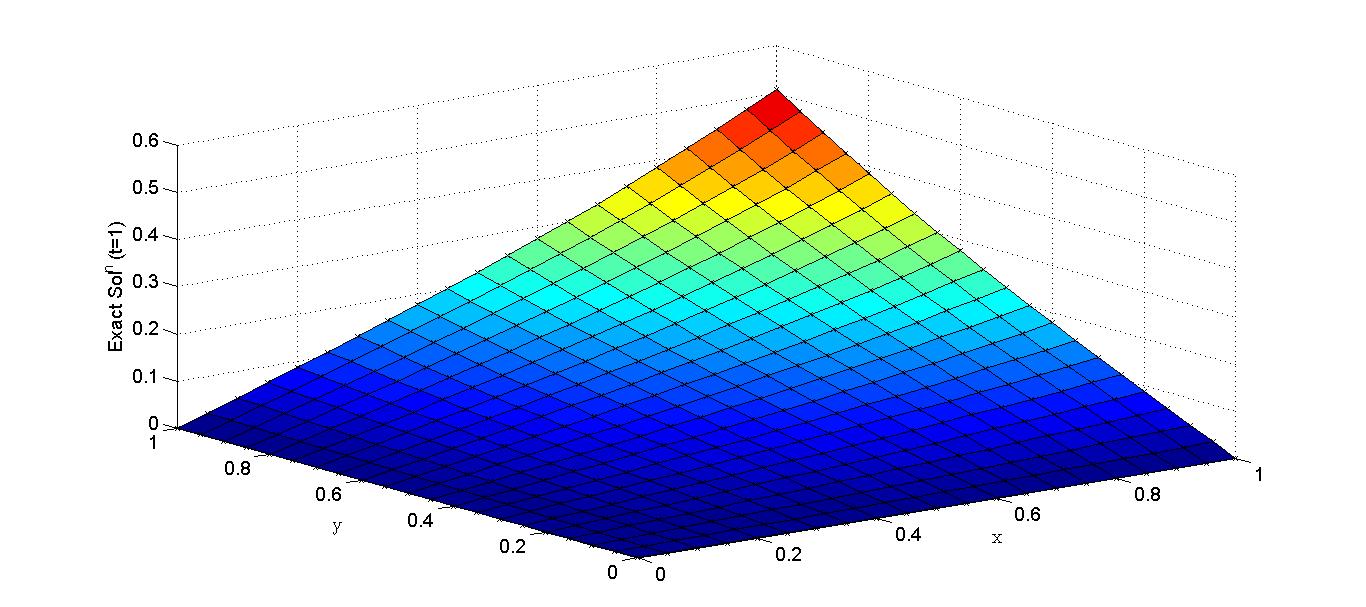}
\includegraphics[height=7.0cm,width=7.6cm]{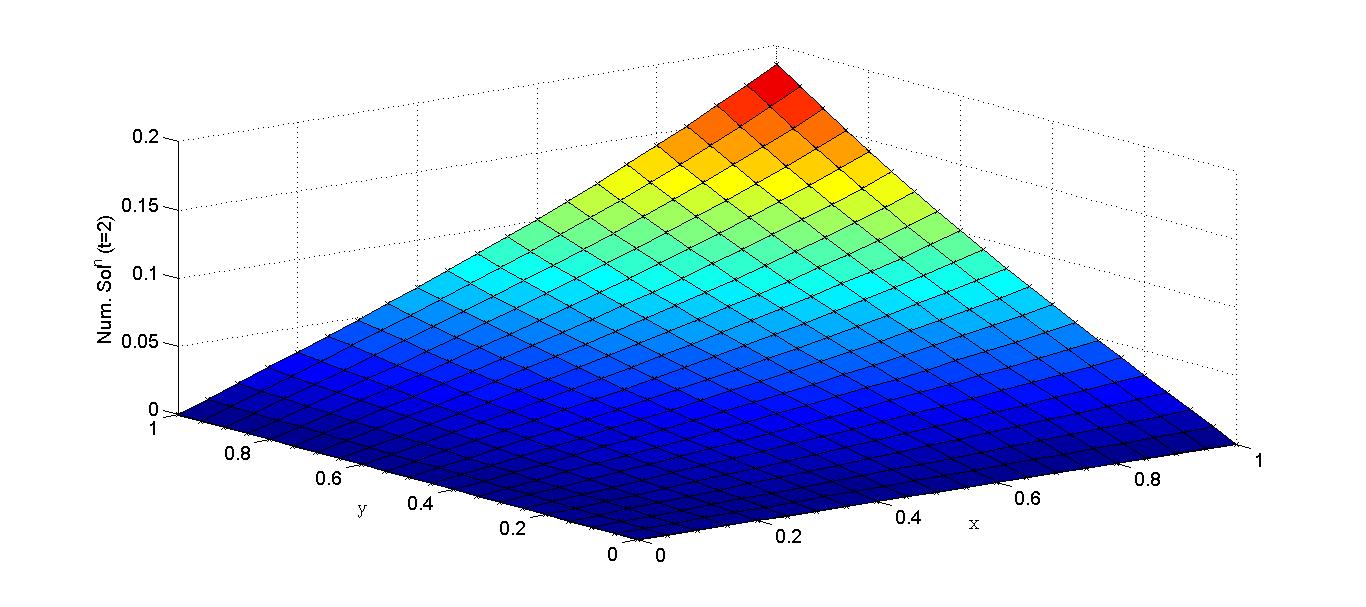}
\includegraphics[height=7.0cm,width=7.6cm]{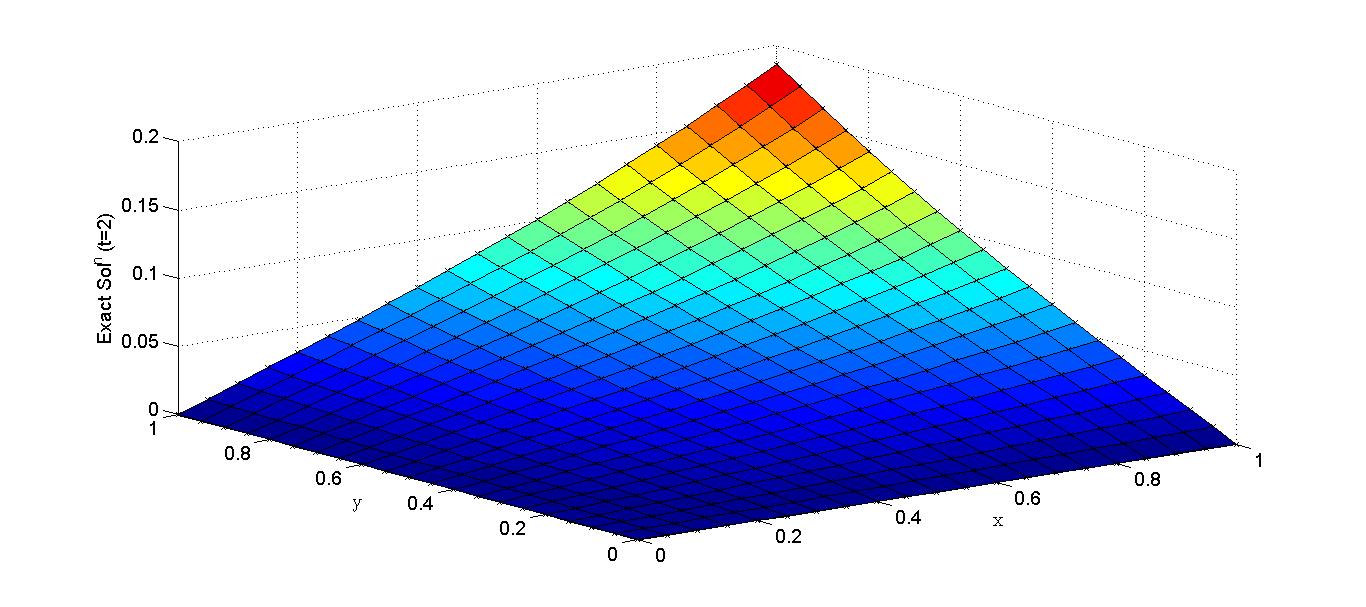}
\includegraphics[height=7.0cm,width=7.6cm]{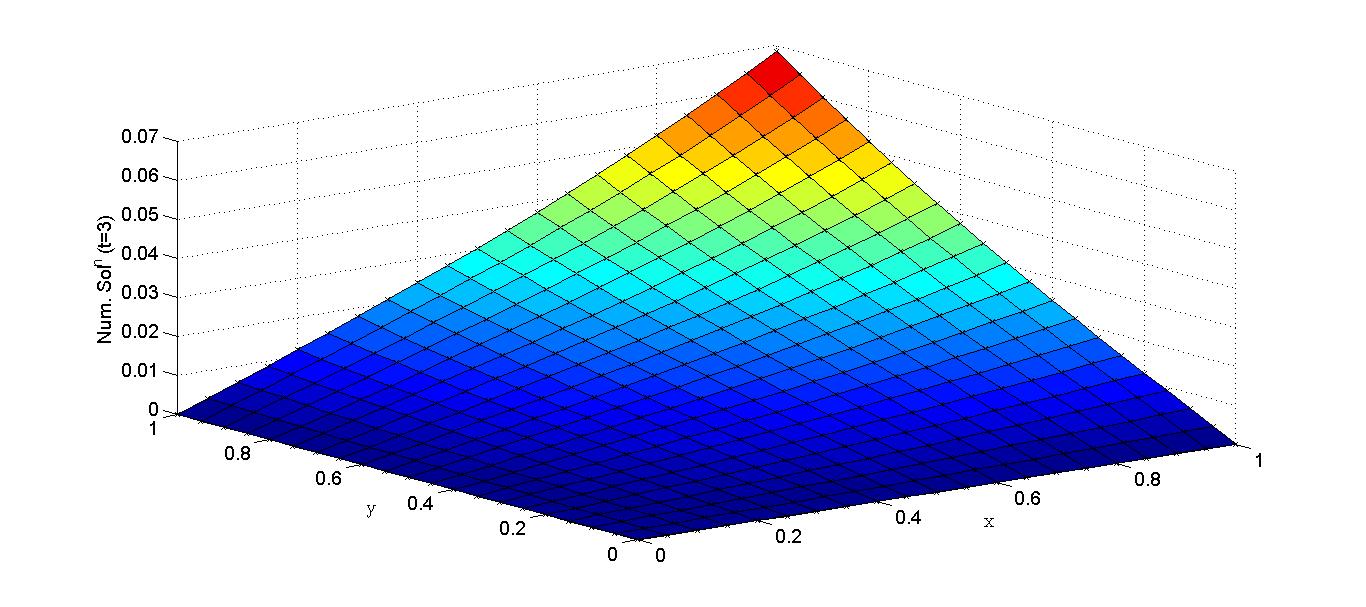}
\includegraphics[height=7.0cm,width=7.6cm]{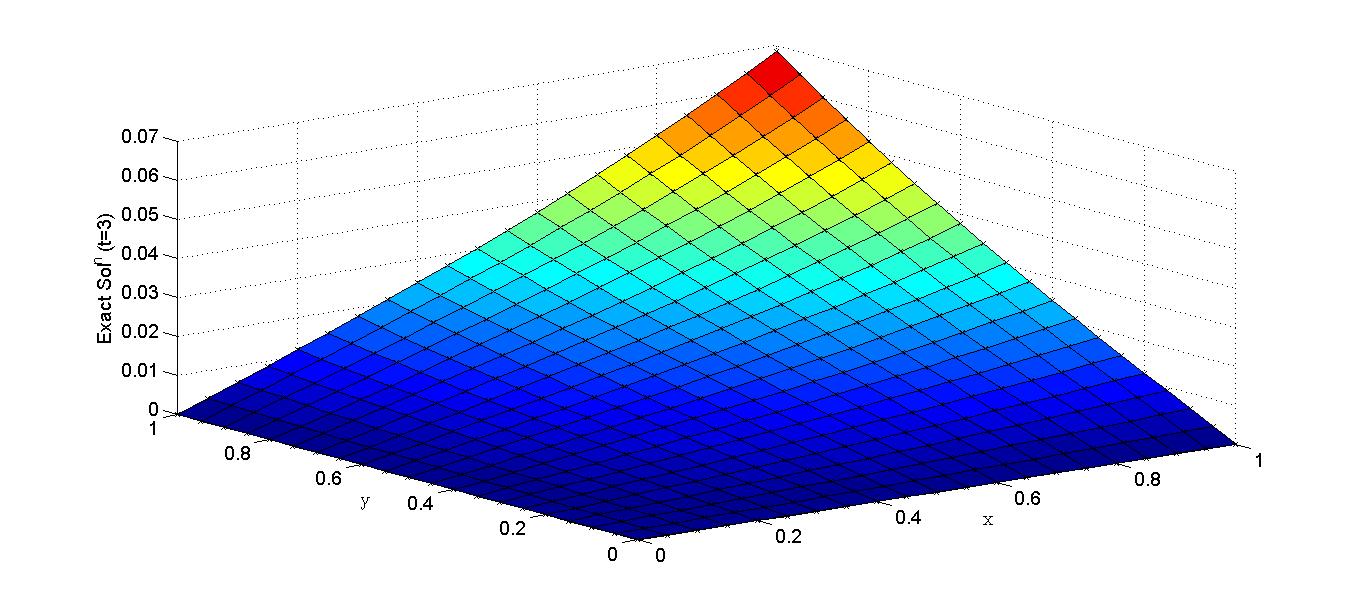}
\caption{Plots of numerical and exact solution at different time levels for Example \ref{ex2}}\label{fig2.1}
\end{figure}

\begin{example}\label{ex2}
Consider the telegraph equation  \eqref{int-tel-eqn1} with $f(x,y,t)=(-2\alpha+\beta^2-1)\exp(-t)\sinh x \sinh y$,
 $\phi(x,y)=\sinh x \sinh y, \psi(x,y)=-\sinh x \sinh y,$ in $\Omega$;$\phi_1(y,t)=0; \phi_2(y,t)=\exp{(-t)}\sinh(1)\sinh y $ for $0\le y\le1$ and $
 \psi_3(x, t)=0; \psi_4(x,t)=\exp {(-t)}\sinh x \sinh (1)$ for $ 0\le x\le 1$.

The exact solution \cite{JPM12} is given by
\begin{equation}\label{ex2-exact-soln}
 u(x,y,t)=\exp{(-t)}\sinh x \sinh y
\end{equation}

The solutions are computed for $\alpha=10,\beta=5$ and $\alpha=10,\beta=0$ with parameters $\Delta t =0.01, 0.001$, $h=0.1, 0.05$ and $p=1$. The computed $L_2,L_\infty$  error norms and CPU time for different time levels $t\le 10$ are compared with the error norms by Mittal and Bhatia \cite{MB14} in Table \ref{tab-ex2.1}, for $\Delta t=0.01$ and $h=0.1$. In Table \ref{tab-ex2.2}, the computed results for $\Delta t=0.001$ and $h=0.05$ are compared with Mittal and Bhatia \cite{MB14} and Jiwari et al. \cite{JPM12}. The findings from the above tables confirms that the proposed results are better than \cite{MB14,JPM12}. The CPU time is slightly more than \cite{MB14} due to selection of SSP-RK54 scheme instead of SSP-RK43 scheme, for time integration. The surface plots of numerical and exact solutions at$ t=1,2,3$ with $\delta t=0.001$ and $h=.05$ are depicted in Fig. \ref{fig2.1}.
\end{example}

\begin{figure}
\centering
\includegraphics[height=7.0cm,width=7.6cm]{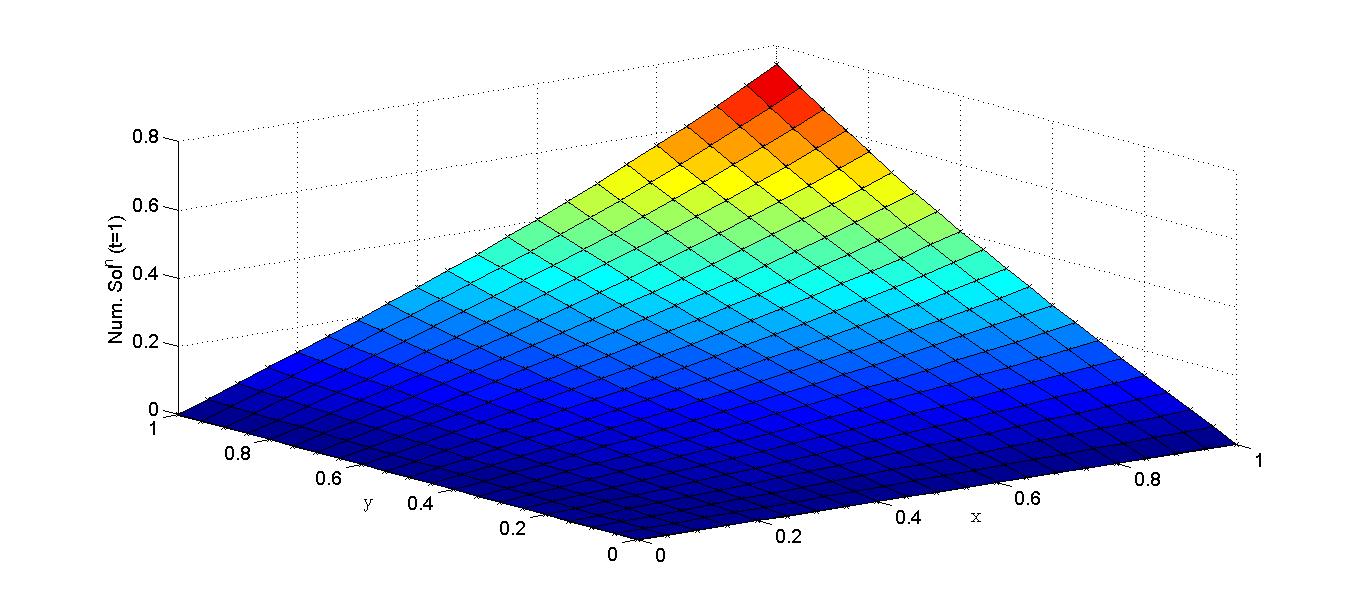}
\includegraphics[height=7.0cm,width=7.6cm]{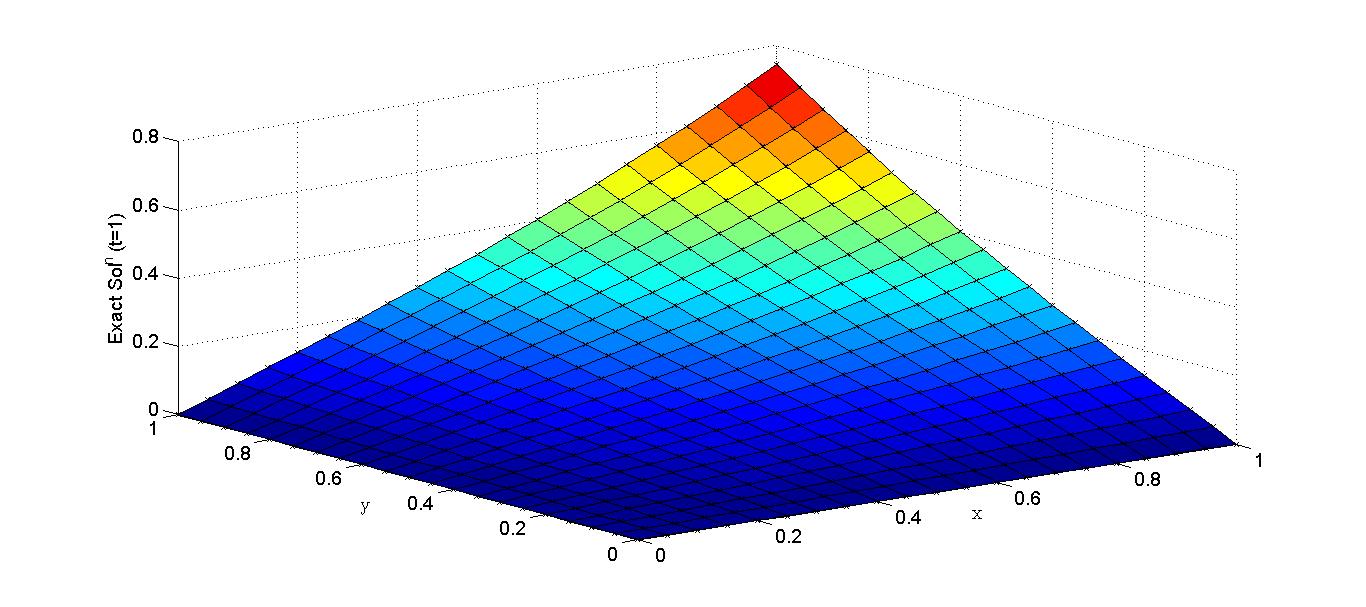}
\includegraphics[height=7.0cm,width=7.6cm]{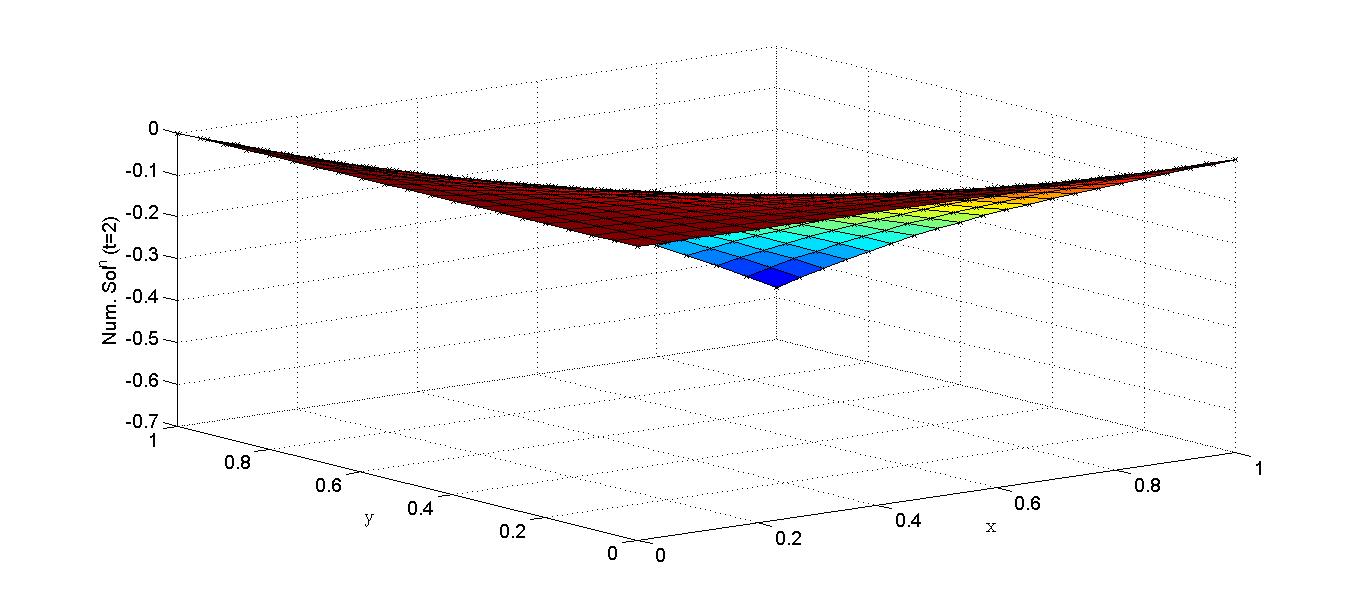}
\includegraphics[height=7.0cm,width=7.6cm]{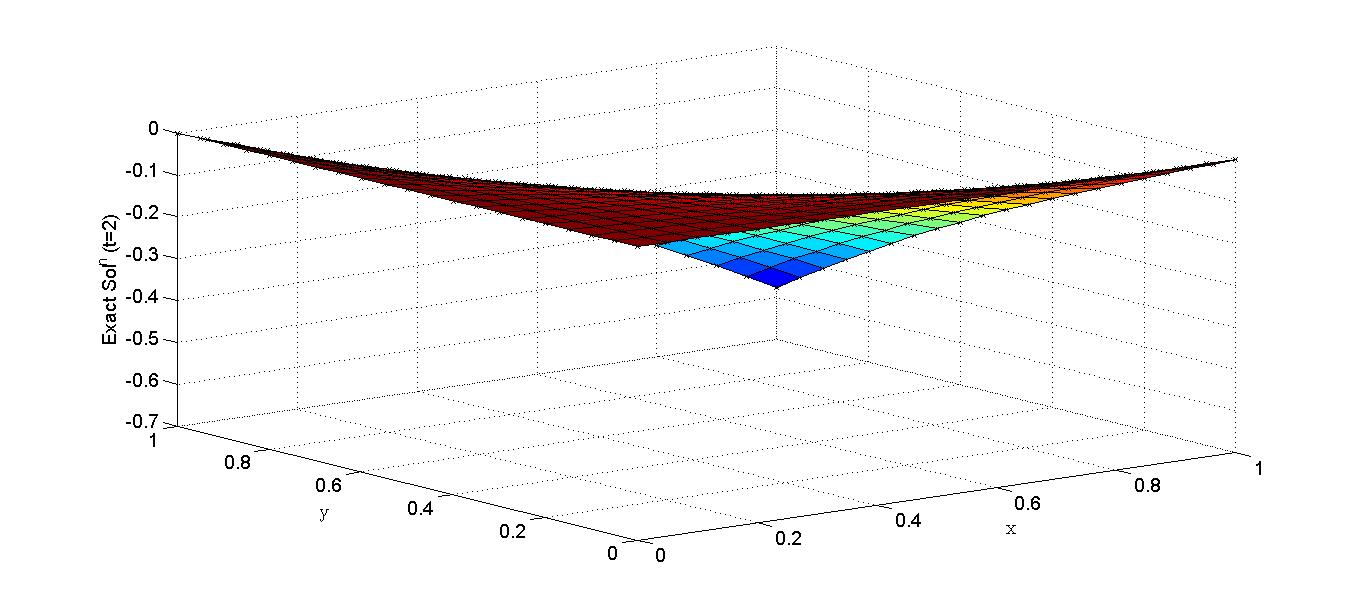}
\includegraphics[height=7.0cm,width=7.6cm]{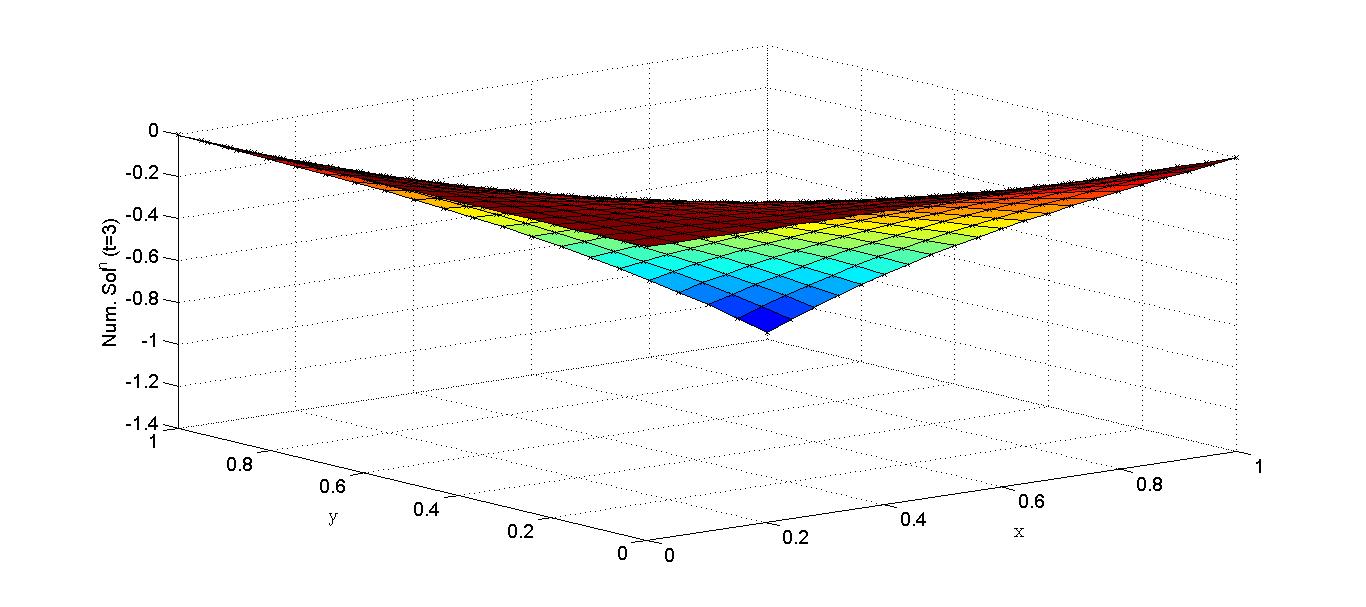}
\includegraphics[height=7.0cm,width=7.6cm]{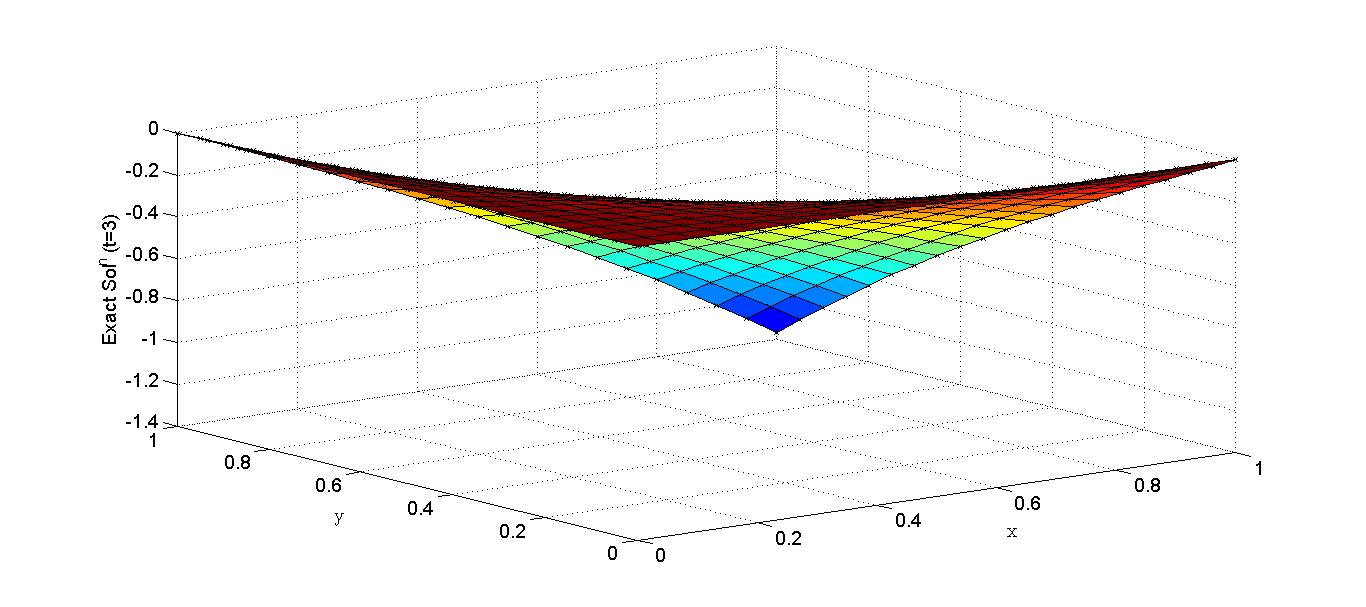}
\caption{Plots of numerical and exact solution at different time levels for Example \ref{ex4}}\label{fig3.1}
\end{figure}

\begin{example}\label{ex4}
 Consider the telegraph equation \eqref{int-tel-eqn1} in the region $\Omega$ with
$f(x,y,t)=(-3 \cos t+2\alpha \sin t +\beta^2 \cos t)\sinh x \sinh y$, and
$ \phi(x,y)= \sinh x \sinh y, \psi(x,y)=0$ in $\Omega$, and $ \phi_1(y,t)=o, \phi_2(y,t)=\cos t \sinh (1)\sinh y $ for $0\le y\le 1$, and $\phi_3(x,t)=0, \psi_4(x,t)=\cos t \sinh x\sinh (1)$ for $0\le x\le 1.$

The exact solution \cite{JPM12} is given by
\begin{equation}\label{ex4-exact-soln}
 u(x,y,t)=\cos t \sinh x \sinh y
\end{equation}
The solution is computed with the parameters $\alpha=10, \beta=5 $ and $\alpha=50, \beta=5$ for the time step $\Delta t=0.001$ and  $h=0.05, p=0.15, 1.$
The computed $L_2, L_\infty$ errors norms and CPU time at different time levels are reported in Table \ref{tab-ex4.1}. It evident that our results are comparably better than the results by Bhatiya and Mittal\cite{MB14}. A comparison of numerical solution with exact solution for $t=1, 2, 3$ is depicted in Fig. \ref{fig3.1}.
\end{example}

\begin{figure}
\centering
\includegraphics[height=7.0cm,width=7.6cm]{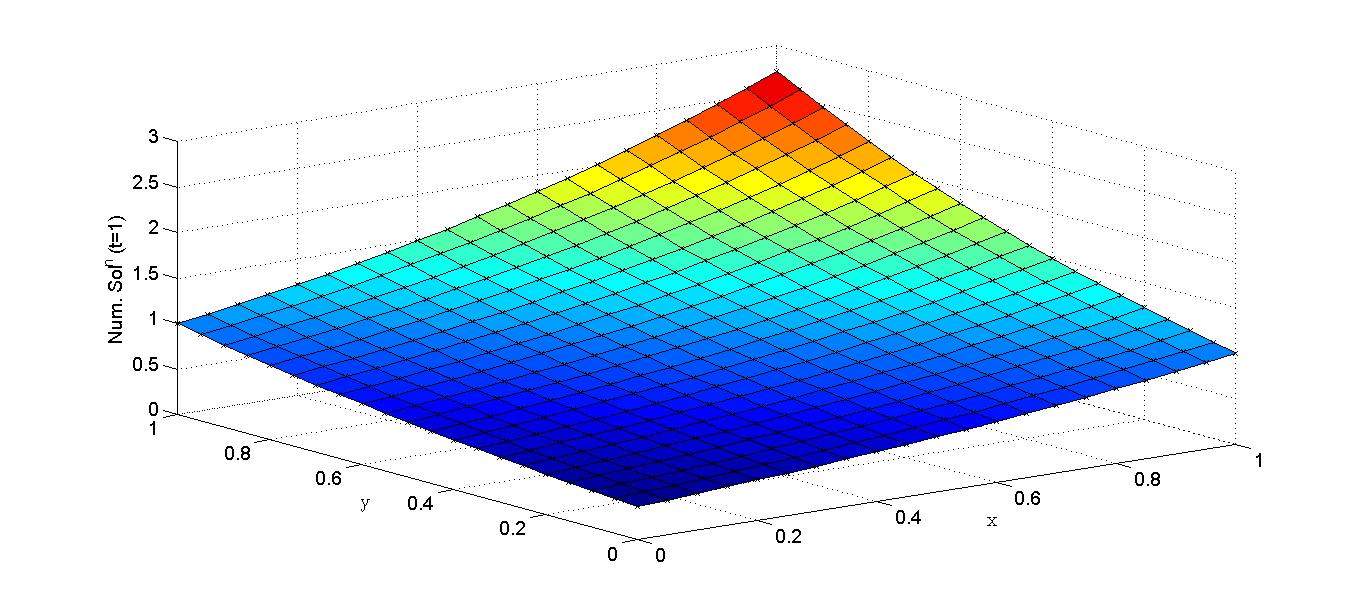}
\includegraphics[height=7.0cm,width=7.6cm]{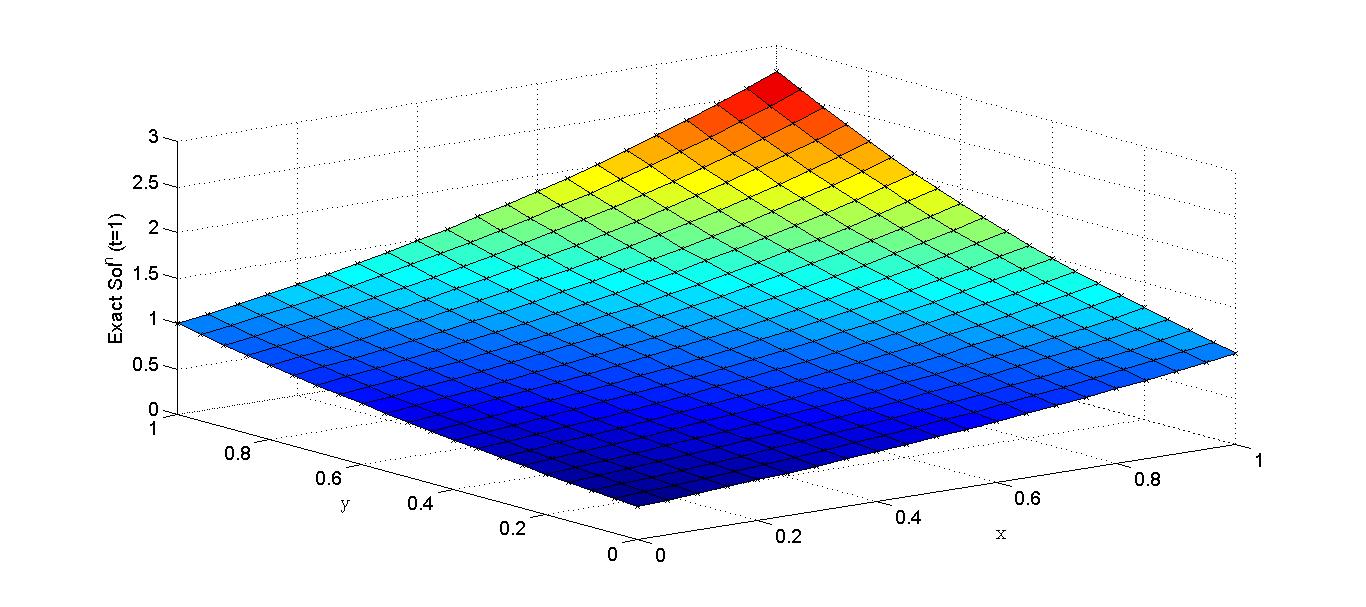}
\includegraphics[height=7.0cm,width=7.6cm]{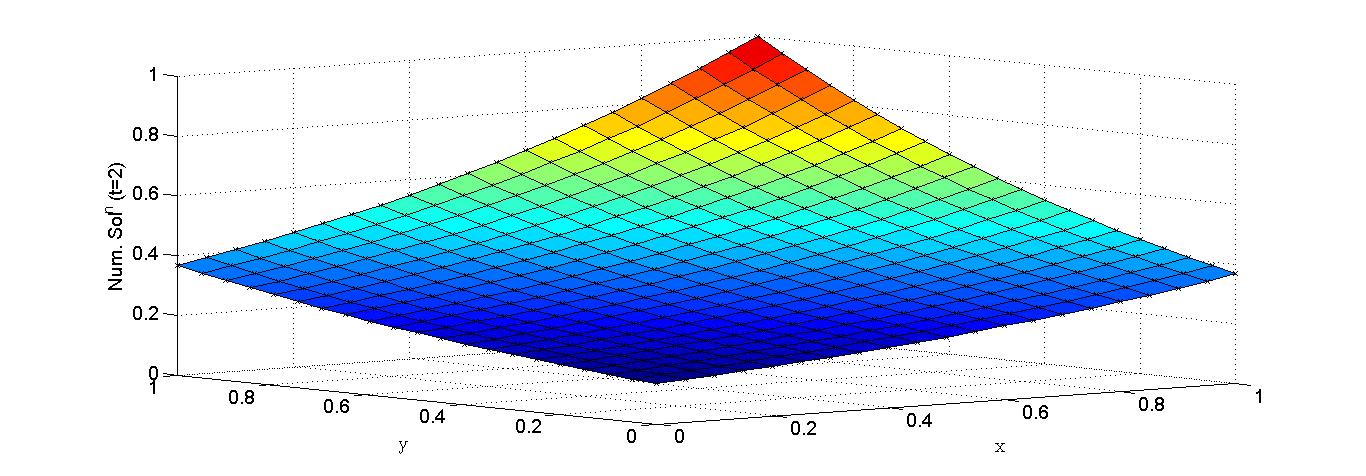}
\includegraphics[height=7.0cm,width=7.6cm]{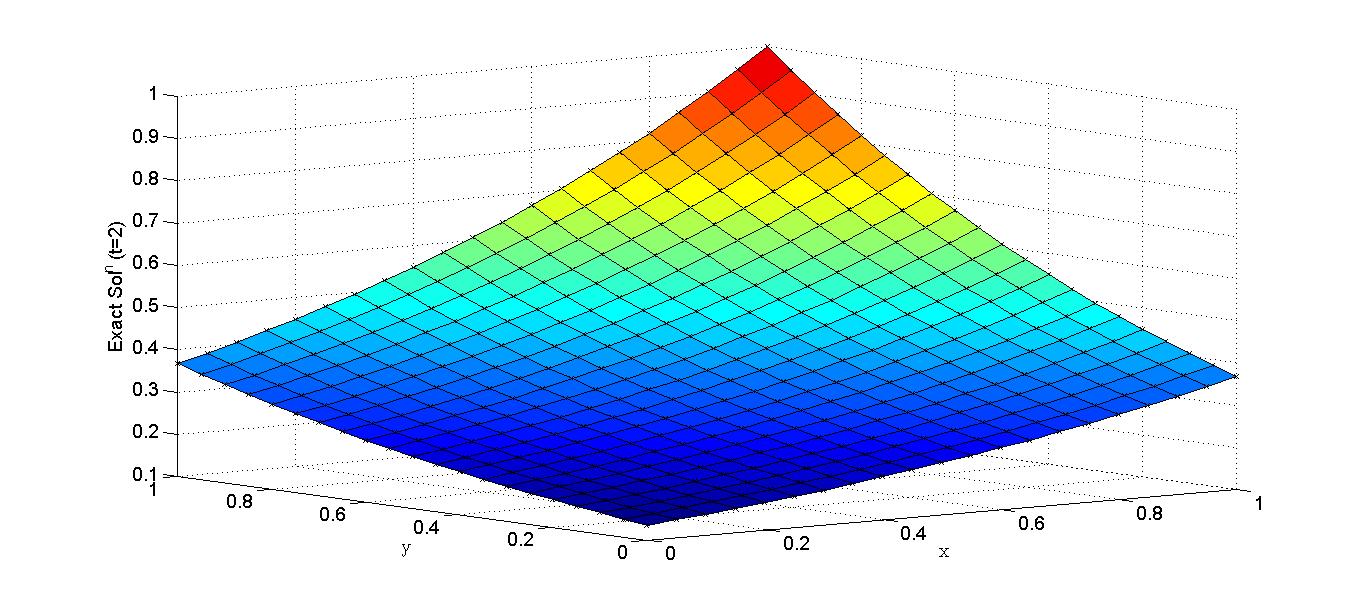}
\includegraphics[height=7.0cm,width=7.6cm]{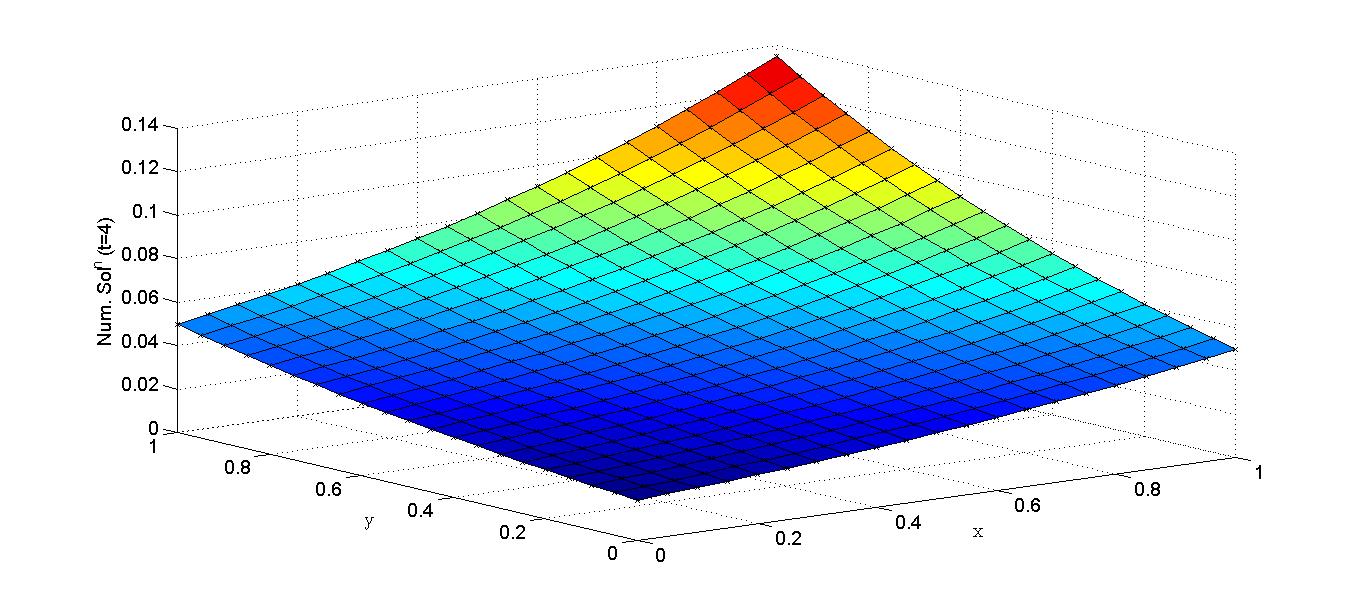}
\includegraphics[height=7.0cm,width=7.6cm]{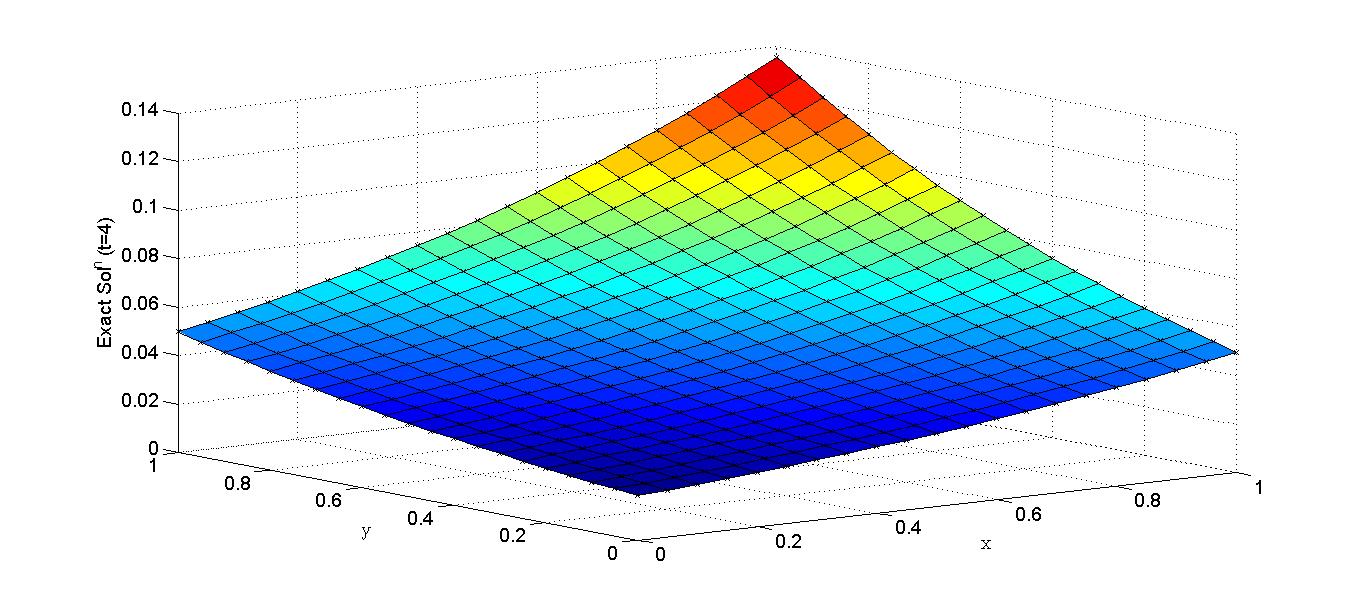}
\caption{Plots of mECDQ solutions and exact solutions at different time levels for Example \ref{ex5}}\label{fig4.1}
\end{figure}

\begin{example}\label{ex5}
Consider the telegraph equation \eqref{int-tel-eqn1} in the region $\Omega$ with $\alpha=1,\beta=1, f(x,y,t)=-2\exp{(x+y-t)}$,
and $ \phi(x,y)= \exp{(x+y)},\psi(x,y)=-\exp{(x+y)}$ in $\Omega$, and the mixed boundary conditions $
 \phi_1(y,t)=\exp {(y-t)}, \phi_2(y,t)=\exp {(1+y-t)} $ for $0\le y\le1 $ and $\psi_3(x, t)=\exp {(x-t)},
\phi_4(x,t)=\exp {(1+x-t)}$ for $ 0\le x\le 1.$
The exact solution \cite{DG10} is given by
\begin{equation}\label{ex5-exact-soln}
 u(x,y,t)=\exp{(x+y-t)}
\end{equation}
The computed results and CPU time are compared with the results by Mittal and Bhatia \cite{MB14} for different space step size $h=0.1, 0.05$ and time step size $\Delta t =0.01, 0.001$, $p=1$ and reported in Table \ref{tab-5.1} and Table \ref{tab-5.2}. The surface plots of the mExp-DQM  solutions and the exact solutions at different time levels $t=1, 2, 4$ is depicted in Fig. \ref{fig4.1}. It is evident accuracy of the proposed results is much better than results of Mittal and  Bhatia \cite{MB14}.
\end{example}

\begin{figure}
\centering
\includegraphics[height=7.0cm,width=7.6cm]{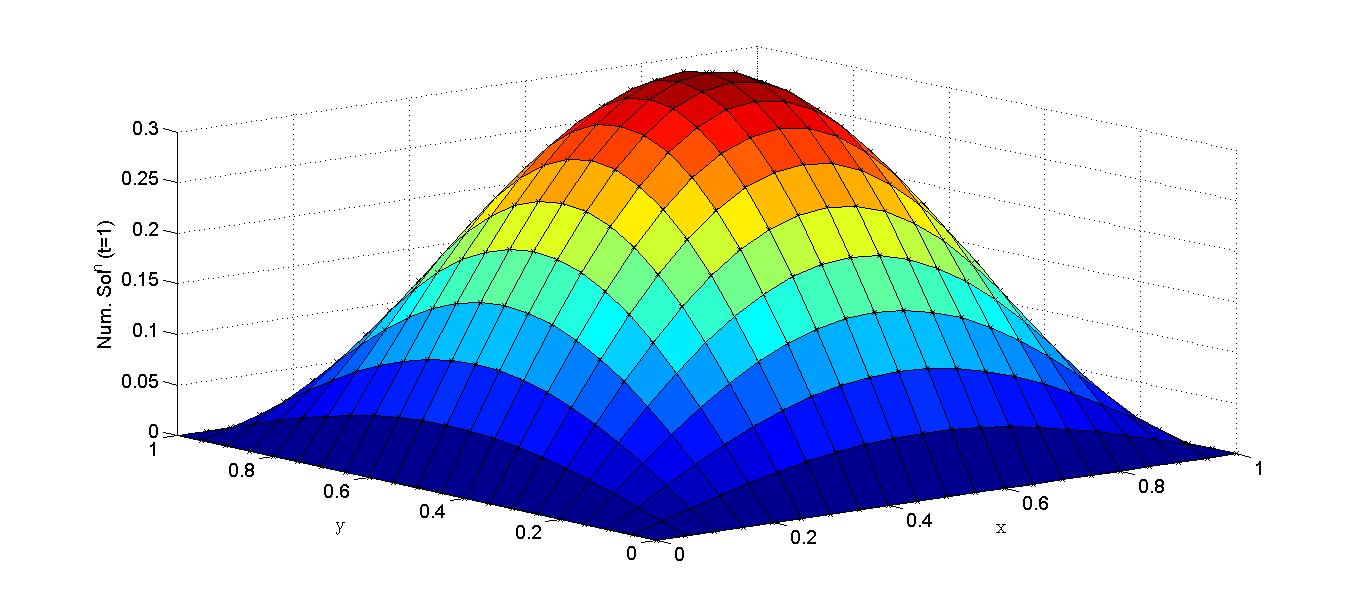}
\includegraphics[height=7.0cm,width=7.6cm]{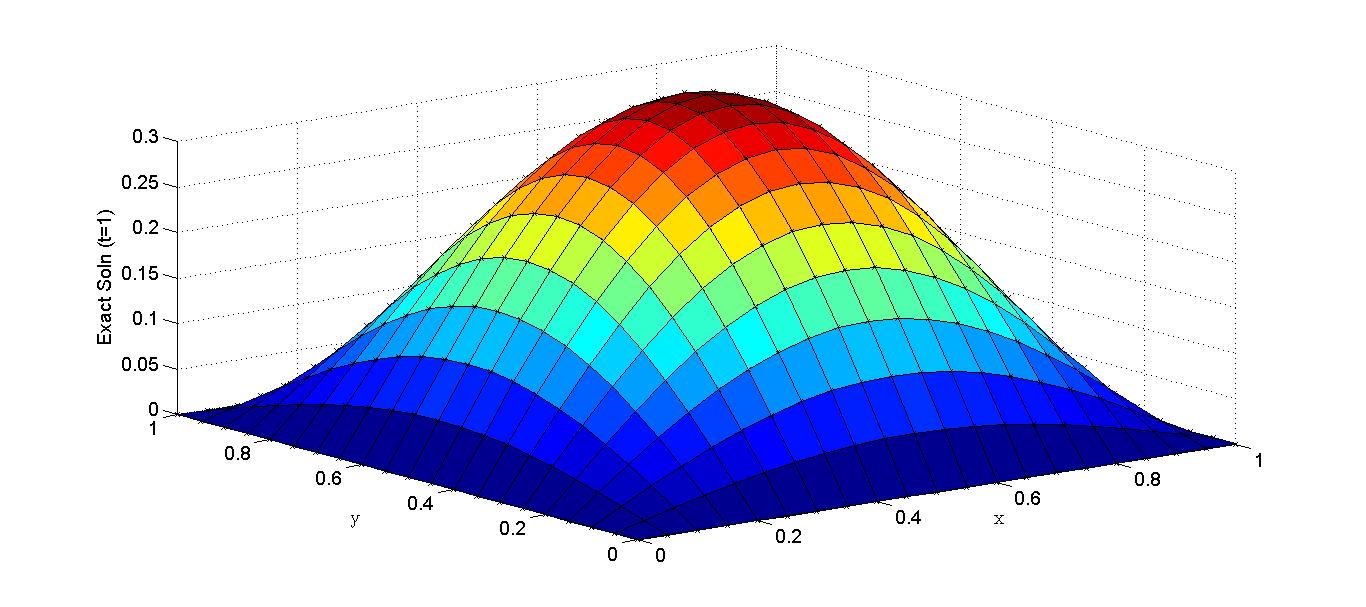}
\includegraphics[height=7.0cm,width=7.6cm]{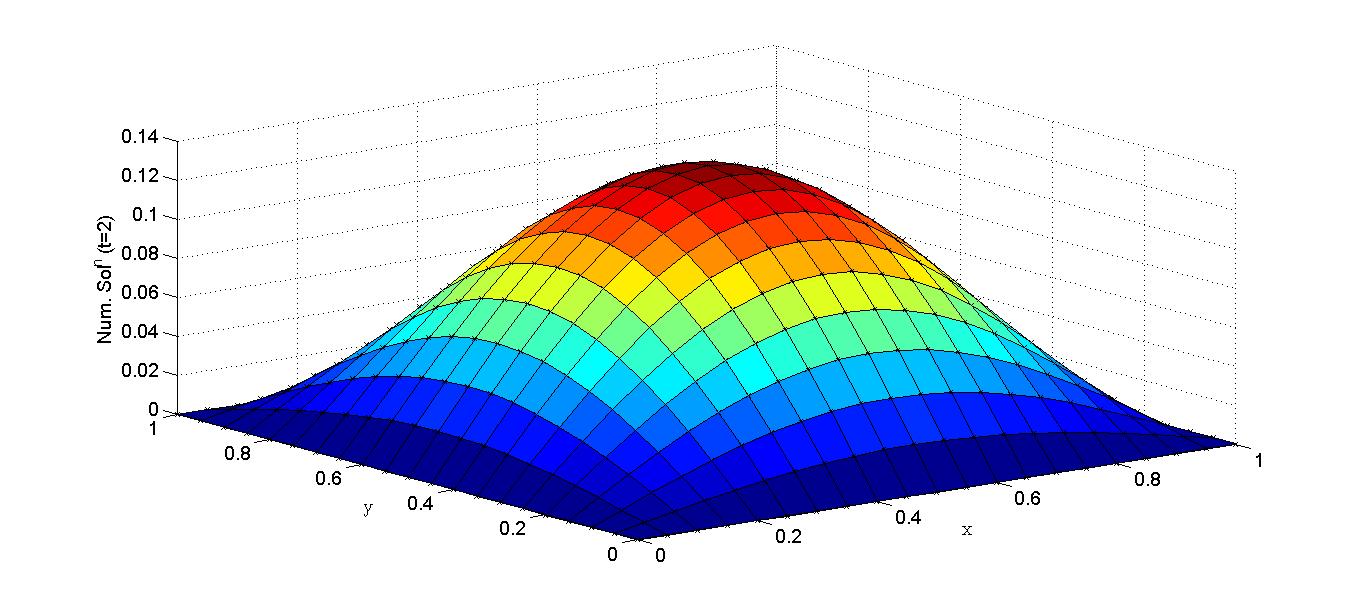}
\includegraphics[height=7.0cm,width=7.6cm]{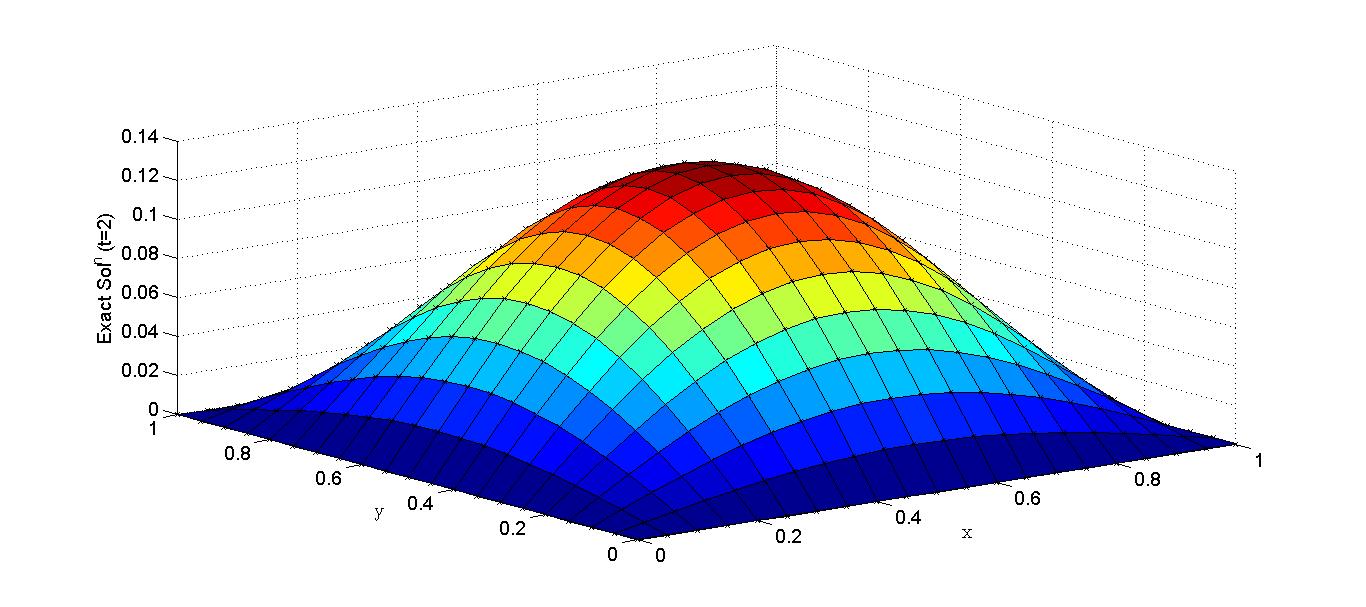}
\includegraphics[height=7.0cm,width=7.6cm]{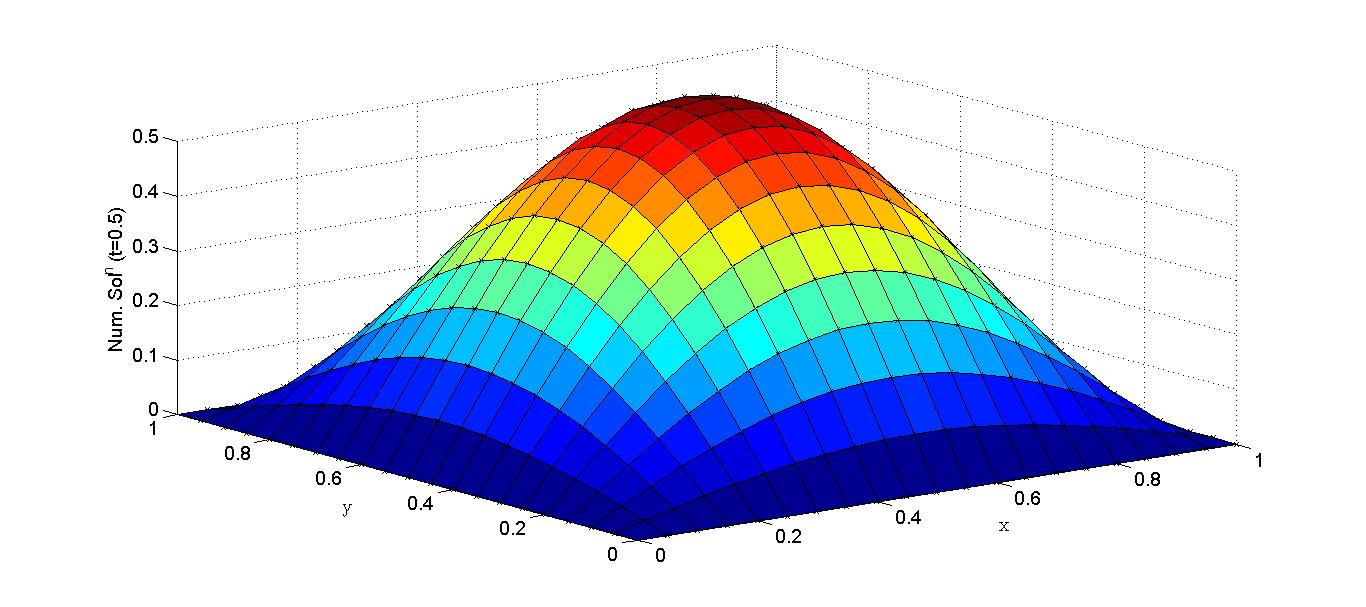}
\includegraphics[height=7.0cm,width=7.6cm]{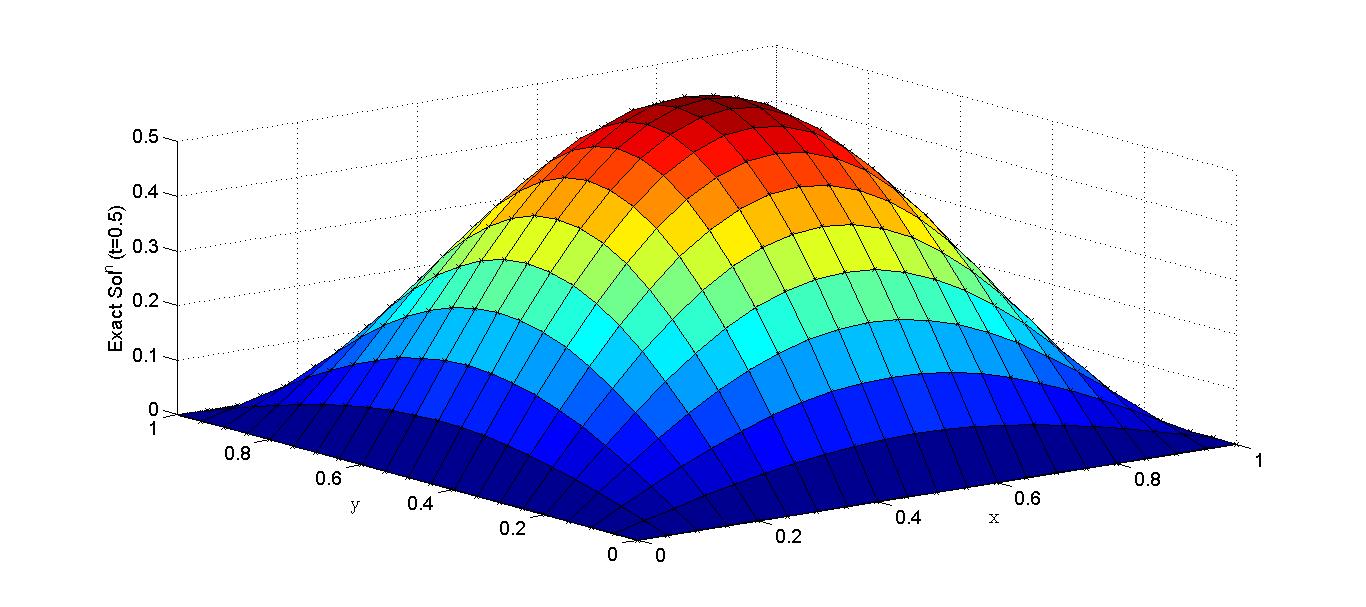}
\caption{Plots of numerical and exact solution at different time levels for Example \ref{ex6}} \label{fig6.1}
\end{figure}

\begin{example}\label{ex6}
The telegraph equation \eqref{int-tel-eqn1} with $\alpha=1, \beta=1$, $f(x,y,t)=2\pi^2 \exp(-t)\sin \pi x \sin {\pi y}$  in region $\Omega, t>0 $ is considered together with $ \phi(x,y)= \sin {\pi x}\sin {\pi y}, \psi(x,y)=-\sin {\pi x}\sin {\pi y}$ in $\Omega$, and the mixed boundary conditions
$ \psi_1(y,t)=\pi \exp {(-t)}\sin{(\pi y)}, \phi_2(y,t)=0$ in  $0 \le y\le 1,$ and
 $\phi_3(x,t)=0,  \psi_4(x,t)=-\pi \exp {(-t)}\sin{(\pi y)}$ in $ 0\le x\le 1$
The exact solution as in \cite{DG10} is given by
\begin{equation}\label{ex5-exact-soln}
 u(x,y,t)=\exp{(-t)}\sin{(\pi x)} \sin{(\pi y)}.
\end{equation}

The solutions are computed in terms of $L_2, L_\infty $ error norms, for $h=0.1, \triangle t =0.01$ and $h=0.05, \triangle t =0.001$ with $p= =0.5, 1$ and reported in Table \ref{tab-ex 6.2} and Table \ref{tab-ex6.1}. The surface plots of numerical and exact solutions at different time levels $t=0.5, 1, 2$ are depicted in Fig. \ref{fig6.1}. The above findings confirms that the proposed mExp-DQM solutions are more accurate as compared to the results by Mittal and Bhatia \cite{MB14}.
\end{example}

\begin{figure}
\centering
\includegraphics[height=7.0cm,width=7.65cm]{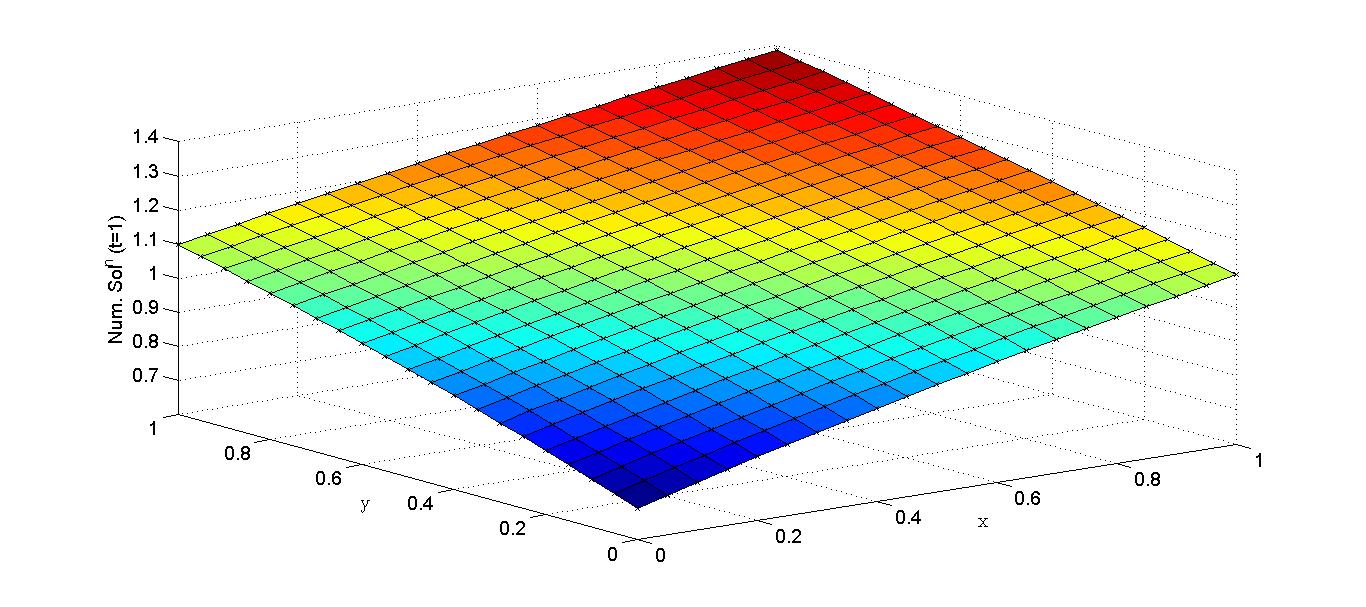}
\includegraphics[height=7.0cm,width=7.65cm]{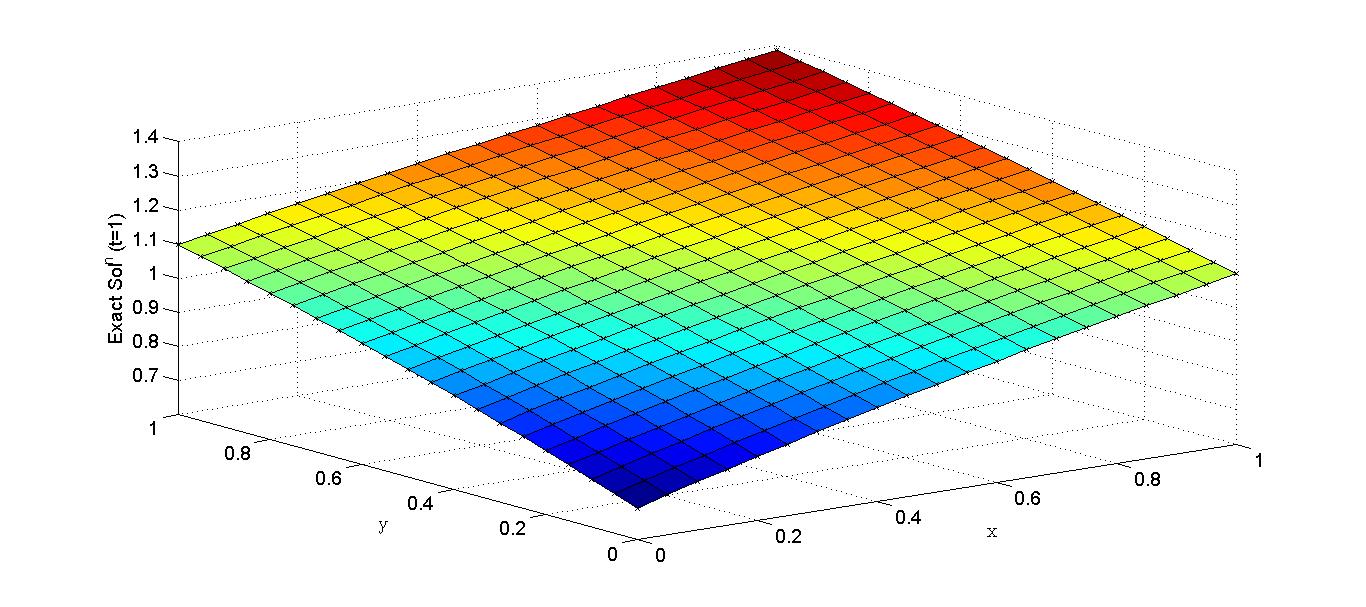}
\includegraphics[height=7.0cm,width=7.65cm]{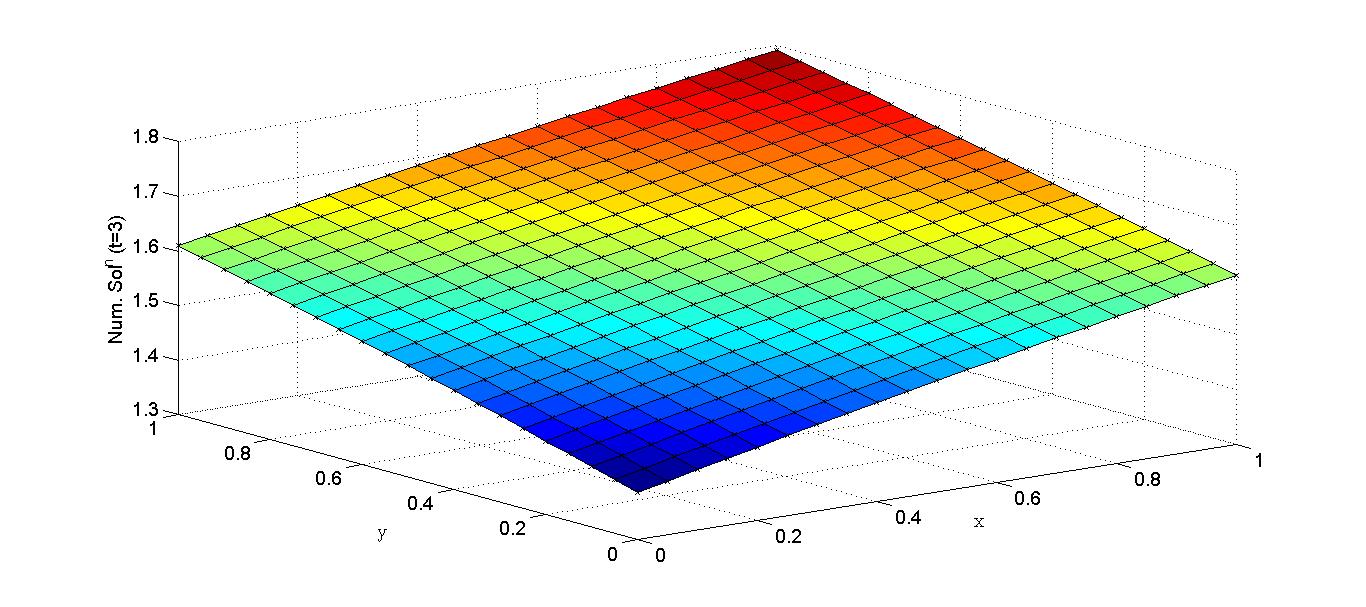}
\includegraphics[height=7.0cm,width=7.65cm]{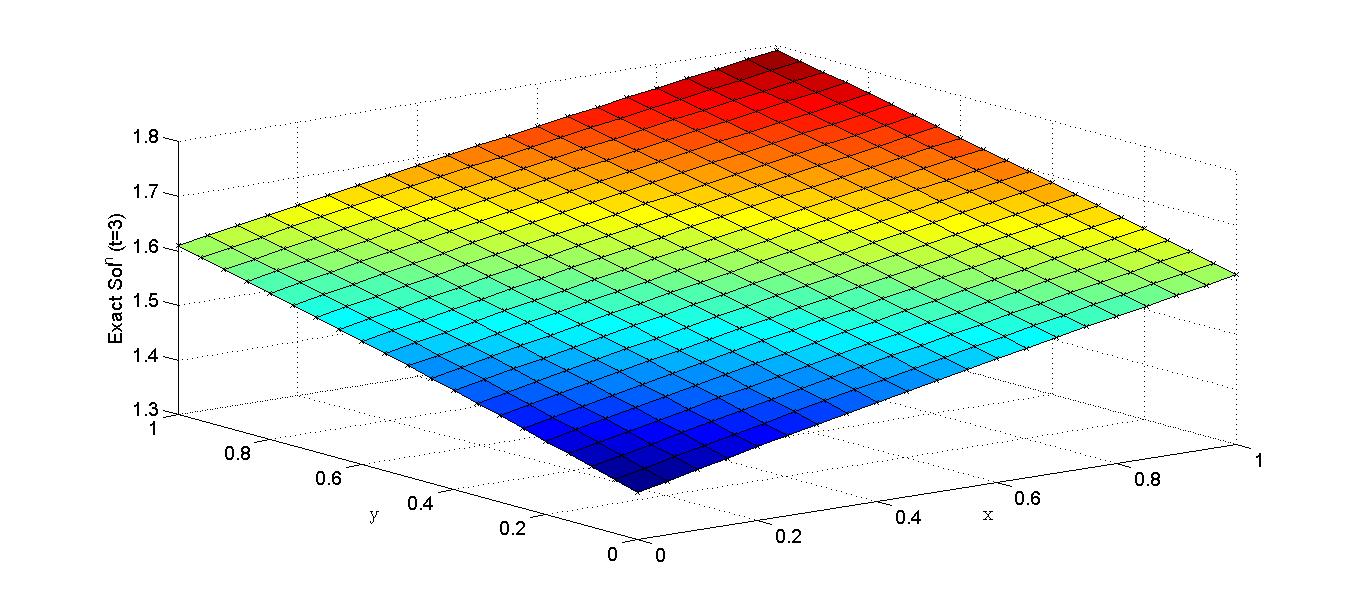}
\includegraphics[height=7.0cm,width=7.65cm]{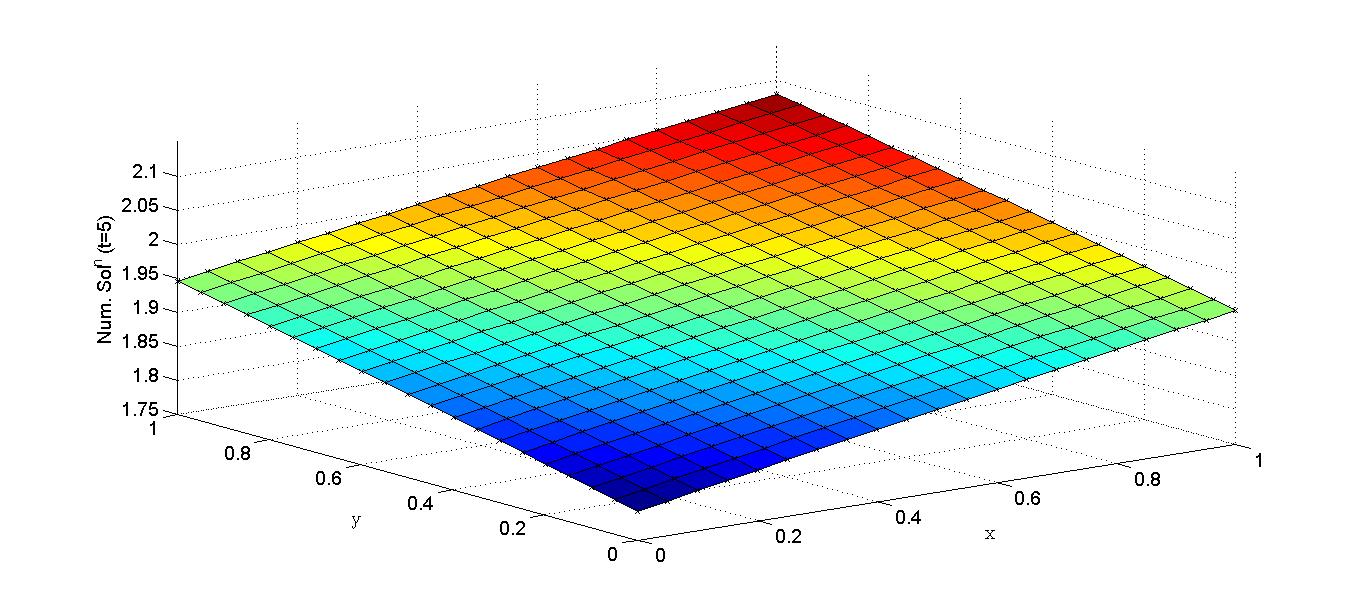}
\includegraphics[height=7.0cm,width=7.65cm]{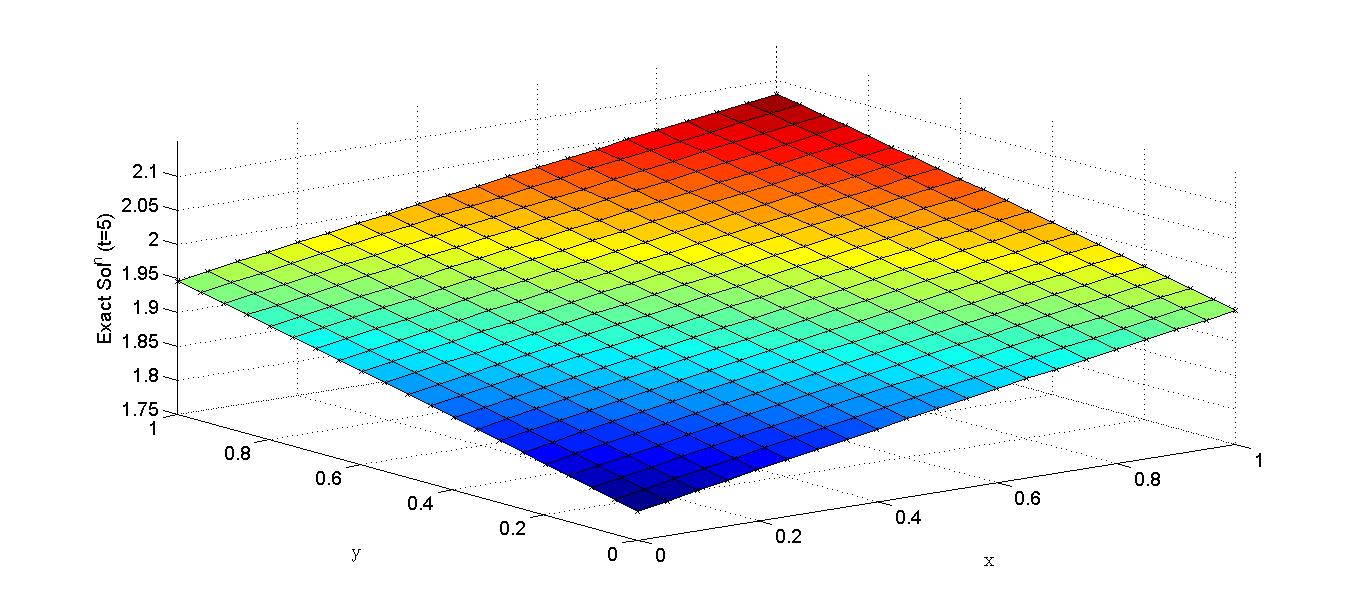}
\caption{Plots of numerical and exact solution at different time levels for Example \ref{ex7}}\label{fig7.1}
\end{figure}

\begin{example}\label{ex7}
 The telegraph equation \eqref{int-tel-eqn1} with  $\alpha=1,\beta=1$ is considered together with
 $\phi(x,y)= \log(1+x+y), \psi(x,y)=\frac{1}{1+x+y}$ in $\Omega$, and the mixed boundary conditions
$ \phi_1(y,t)=\log(1+y+t), \psi_2(y,t)=\frac{1}{2+y+t}$ for  $0\le y\le 1$ and $
\psi_3(x,t)=\frac{1}{1+x+t}, \phi_4(x,t)=log(2+x+t)$ for $0\le x\le 1.$ The exact solution as given in \cite{DG10} is:
\begin{equation}\label{ex7-exact-soln}
 u(x,y,t)=\log(1+x+y+t)
\end{equation}
where the function $f(x,y,t)$ can be extracted from the exact solution.

The solution is computed for $p=1$, $\Delta t=0.001, h=0.05$ in the region $\Omega$ in terms of $L_2, L_\infty$ and relative error norms. The computed results are compared with the results by Mittal and Bhatia \cite{MB14} and Dehghan and Ghesmati \cite{DG10a}, reported in Table \ref{tab-ex7.0}. It is evident from Table  \ref{tab-ex7.0} that the accuracy of mExp-DQM results is much better than the accuracy in the results of \cite{MB14}, and \cite{DG10a} for large $t$. The surface plots of numerical and exact solutions at different time levels $t=1,2,3$ are depicted in Fig. \ref{fig7.1}.

\end{example}

\section{Conclusion} \label{sec-conclu}
In this paper, we have developed a new differential quadrature method based on modified exponential cubic B-splines as a set of basis functions, and so, we called it modified exponential cubic B-spline differential quadrature  method (mExp-DQM). The developed mExp-DQM with SSP-RK54 scheme is implemented for second order hyperbolic telegraph equation in $(2+1)$ dimension subject to initial conditions, and each type of boundary conditions: Drichlet, Neumann, mixed boundary conditions.

The results are compared with the recent results by Mittal and Bhatia \cite{MB14} and Jiwari et al. \cite{JPM12}. It is evident that the accuracy of propped results is good as compared to very recent and accurate results due to \cite{MB14,JPM12}. The CPU time is more than \cite{MB14}, while very less in comparison to \cite{JPM12}. Finally, we conclude that the proposed mExp-DQM results with suitable value of free parameter $p$ produces comparatively good results than \cite{MB14,JPM12}.

\section*{Acknowledgement}
\noindent 
Pramod Kumar would like to thanks BBA University Lucknow, India for financial assistance to carry out the research work.

\newpage

\section{List of Tables and Figures}

\begin{table}[!htbp]
\caption{Comparison of the mExp-DQM solutions of Example \ref{ex1} with $\Delta t=0.01$ and $h=0.1, p=1$}\label{tab-ex1.1}
\vspace{.2cm}
\centering
\begin{tabular}{ll*9l}
\toprule
$t$&\hspace{0.25cm}& \multicolumn{4}{c}{ mExp-DQM }  &{}& \multicolumn{4}{c}{ Mittal Bhatia \cite{MB14} }   \\ \cline{2-5} \cline{7-10}
{} &$L_2$ &$L_\infty$& $R_e$ &$CPU(s)$&{}&$L_2$ &$L_\infty$ &$R_e$&$CPU(s)$\\
\midrule				
1	$\quad$&3.7330E-06	$\quad$&4.5492E-06	$\quad$&2.7069E-04	$\quad$&0.031	$\quad$&&9.9722E-04	$\quad$&2.2746E-03	$\quad$&5.9762E-03	$\quad$&0.08\\
2	&4.4842E-06	&5.6294E-06	&4.2217E-04	&0.062	&&1.0926E-03	&2.8706E-03	&8.5019E-03	&0.11\\
3	&3.7742E-06	&6.3374E-06	&1.4916E-04	&0.109	&&2.2877E-04	&6.0818E-04	&7.4720E-04	&0.14\\
5	&4.4186E-06	&5.3912E-06	&5.9036E-04	&0.203	&&1.1562E-03	&2.9942E-03	&1.2767E-03	&0.20\\
7	&3.2109E-06	&3.7239E-06	&1.6834E-04	&0.312	&&7.2867E-04	&1.8781E-03	&3.1572E-03	&0.26\\
10	&3.1806E-06	&3.7506E-06	&1.4949E-04	&0.374	&&5.8889E-04	&1.5158E-03	&2.2874E-03	&0.34\\

\bottomrule
\end{tabular}
\end{table}
 \begin{table}[!htbp]
\caption{Comparison of the  mExp-DQM solutions of Example \ref{ex1} with $\Delta t=0.001$; $ h=0.05$ and $p=1,0.15$}\label{tab-ex1.2}
\vspace{.2cm}
\centering
\begin{tabular}{*4lcl*3lc*4lc}
\toprule
$t$&   \multicolumn{4}{l}{  mExp-DQM (p=1)}  &{}& \multicolumn{4}{l}{  mExp-DQM (p=0.15)} &{}&\multicolumn{4}{l}{  Mittal and Bhatia \cite{MB14} }   \\ \cline{2-5} \cline{7-10} \cline{12-15}
{} &$L_2$ &$L_\infty$ &$R_e$&CPU(s)&{}&$L_2$ &$L_\infty$ &$R_e$&CPU (s)&{}&$L_2$ &$L_\infty$ &$R_e$&CPU(s)\\
\midrule			
1	~~&3.5715E-07	~&5.8718E-07	~~&5.0162E-05	&2.26	&&3.5729E-07	~~&5.8736E-07	~~&5.0182E-05	&2.26	&&9.8870E-05~&2.4964E-04	 ~&6.2977E-04	 &0.78\\
2	&4.4969E-07	&6.7211E-07	&8.1823E-05	&4.52	&&4.4977E-07	&6.7225E-07	&8.1838E-05	&4.52	&&1.2148E-04	&3.2296E-04	&1.0025E-03	&1.30\\
3	&7.8128E-07	&1.2228E-06	&5.9879E-05	&6.79	&&7.8151E-07	&1.2231E-06	&5.9897E-05	&6.79	&&3.7627E-05	&9.9310E-05	&1.3078E-04	&1.70\\
5	&3.6743E-07	&4.4790E-07	&9.8297E-05	&11.17	&&3.6749E-07	&4.4797E-07	&9.8311E-05	&11.17	&&1.2762E-04	&3.3205E-04	&1.5411E-03	&3.00\\
7	&5.0400E-07	&8.8032E-07	&5.0732E-05	&15.81	&&5.0417E-07	&8.8056E-07	&5.0749E-05	&15.81	&&6.7672E-05	&1.7679E-04	&3.0892E-04	&3.30\\
10	&5.7992E-07	&9.9151E-07	&5.2483E-05	&22.60	&&5.8011E-07	&9.9178E-07	&5.2500E-05	&22.60	&&5.1764E-05	&1.3521E-04	&2.1245E-04	&5.20\\
\bottomrule
\end{tabular}
\end{table}

\begin{table}[!htbp]
\caption{Comparison of the mExp-DQM solutions of Example \ref{ex2} with $\alpha=10,\beta=5,\Delta t=0.01,h=0.1 $ and $p=1$}\label{tab-ex2.1}
\vspace{.2cm}
\centering
\begin{tabular}{*9lr}
\toprule
$t$&  \multicolumn{4}{l}{Mittal and Bhatia \cite{MB14} }  &{}& \multicolumn{4}{l}{mECDQ method}   \\ \cline{2-5} \cline{7-10}
{}&$L_2$ &$L_\infty$ & $R_e$&$CPU(s)$&{}&$L_2$ &$L_\infty$ & $R_e$&$CPU(s)$\\
\midrule
0.5 $\quad$& 8.3931E-04	$\quad$&3.3019E-03	$\quad$&2.8902E-03	$\quad$&0.13	&&6.8998E-06	$\quad$&1.0168E-05	$\quad$&2.8749E-04	$\quad$&0.016\\
1	& 6.0254E-04	&2.0597E-03	&3.4208E-03	&0.16	&&5.3522E-06	&7.0133E-06	&3.6767E-04	&0.046\\
2	& 2.4167E-04	&7.6531E-04	&3.7297E-03	&0.19	&&2.2337E-06	&2.8534E-06	&4.1711E-04	&0.078\\
3	& 8.9534E-05	&2.7920E-04	&3.7937E-03	&0.24	&&8.3375E-07	&1.0585E-06	&4.2747E-04	&0.141\\
5	& 1.2168E-05	&3.7800E-05	&3.8097E-03	&0.34	&&1.1352E-07	&1.4389E-07	&4.3005E-04	&0.218\\
\bottomrule
\end{tabular}
\end{table}

 \begin{table}[!htbp]
\caption{Comparison of mECDQ solutions of Example \ref{ex2} with $\delta t=0.001,\alpha=10, \beta=0, 5, p=1$ and $h=0.05$}\label{tab-ex2.2}
\vspace{.2cm}
\centering
\begin{tabular}{l*3lc*4lc*2lc}
\toprule
$t$&  \multicolumn{4}{l}{ mECDQ method }  &{}& \multicolumn{4}{l}{Mittal and Bhatia \cite{MB14} } &{}& \multicolumn{2}{l}{Jiwari et al. \cite{JPM12} }  \\ \cline{2-5} \cline{7-10} \cline{12-13}
$\beta=5$ &$L_2$&	$L_\infty$&	$R_e$&	CPU(s) & &	$L_2$&	$L_\infty$	 &$R_e$&	CPU(s)& &	$R_e$ &	CPU(s)\\
\midrule																															
0.5		~~&8.1273E-07	~~~~&1.3152E-06	~~~~&6.6847E-05	&1.279	&&1.0690E-04	~~~~&2.4738E-04	~~~~&1.1088E-04	~~&0.47	&&1.1185E-04~~&	6\\
1		&5.8429E-07	&8.3976E-07	&7.9233E-05	&2.496	&&1.5293E-05	&3.3082E-04	&1.3266E-04	&1.10	&&1.8051E-04&	12\\
2		&2.3507E-07	&3.2200E-07	&8.6737E-05	&5.896	&&4.6468E-05	&1.1380E-05	&3.1954E-04	&1.10	&&4.7289E-04&	25\\
3		&8.8032E-08	&1.1937E-07	&8.8297E-05	&7.488	&&2.1994E-05	&4.3577E-05	&1.3024E-04	&2.80	&&1.2656E-04&	37\\
5		&1.1979E-08	&1.6202E-08	&8.8694E-05	&12.402	&&2.7151E-06	&5.4141E-06	&1.4439E-04	&4.30	&&9.2770E-04&	62\\
$\beta=0$& &&&&&&&&&&& \\											
0.5		&7.2898E-07	&9.0815E-07	&5.9958E-05	&1.279	&&9.2959E-05	&4.2348E-04	&3.4675E-04	&0.52	&&1.1198E-04&	6\\
1		&8.0739E-07	&1.0270E-06	&1.0949E-04	&2.511	&&6.3652E-05	&2.5838E-04	&3.9146E-04	&0.98	&&1.8635E-04&	12\\
2		&5.7525E-07	&7.2622E-07	&2.1226E-04	&5.007	&&2.5540E-05	&9.5843E-05	&4.2739E-04	&1.80	&&5.1797E-04&	25\\
3		&3.1155E-07	&3.9340E-07	&3.1248E-04	&7.394	&&9.9234E-06	&3.5340E-05	&4.5140E-04	&2.20	&&1.4412E-04&	37\\
5		&6.7799E-08	&8.5767E-08	&5.0198E-04	&12.470	&&1.5116E-06	&4.8043E-06	&5.0758E-04	&4.50	&&1.0883E-04&	62\\
\bottomrule
\end{tabular}
\end{table}

 \begin{table}[!htbp]
\caption{Comparison of mECDQ solutions of Example \ref{ex4} with $\delta t=0.001,\alpha=10, 50,\beta=5$ and $h=0.05$}\label{tab-ex4.1}
\vspace{.2cm}
\centering
\begin{tabular}{*3lc*2lr*3lr}
\toprule
$t$&  \multicolumn{6}{l}{ mECDQ method }  &{}& \multicolumn{3}{l}{Mittal and Bhatia \cite{MB14} } \\ \cline{2-7} \cline{8-11}
$\alpha=10$ &$L_2: p=0.015$ &	$L_\infty: p=0.015$&	CPU(s) $~~~~$&	$L_2: p=1$	& $L_\infty: p=1$	& CPU(s) & &	$L_2$	& $L_\infty$ &	CPU(s)\\
\midrule																			
0.5		$~~~~$&2.0862E-06$~~~~$&	2.8531E-06$~~~~$&	1.294	&2.0861E-06$~~~~$&	2.8527E-06$~~~~$&	1.310	&&1.070E-04	$~~~~$&3.756E-04$~~~~$&	0.57\\
1		& 2.5046E-06&	3.2481E-06&	2.62	&2.5045E-06&	3.2479E-06&	2.608	&&1.717E-04	&5.640E-04&	0.92\\
2		& 1.3896E-06&	1.7942E-06&	6.115	&1.3896E-06&	1.7942E-06&	5.179	&&1.647E-04	&5.130E-04&	1.20\\
3		& 1.4008E-06&	2.2256E-06&	7.722	&1.4006E-06&	2.2252E-06&	7.722	&&8.986E-06	&1.956E-05&	2.30\\
5		& 1.6566E-06&	2.1478E-06&	12.885	&1.6566E-06&	2.1478E-06&	12.901	&&1.774E-04	&5.563E-04&	4.10\\
7		& 2.5344E-06&	3.2884E-06&	18.142	&2.5342E-06&	3.2882E-06&	18.008	&&1.420E-04	&4.723E-04&	5.40\\
10		& 1.8983E-06&	2.4641E-06&	25.631	&1.8983E-06&	2.4640E-06&	25.646	&&1.224E-04	&4.122E-04&	7.40\\
$\alpha=50$	& &&&&&&&&& \\
0.5		&2.2128E-06&	3.2835E-06&	1.544	&2.2127E-06&	3.2833E-06&	1.294	&&9.880E-05	&3.696E-04&	0.57\\
1		&3.2434E-06&	4.4892E-06&	2.574	&3.2433E-06&	4.4891E-06&	2.605	&&1.677E-04	&5.687E-04&	0.94\\
2		&2.4069E-06&	3.2575E-06&	5.194	&2.4069E-06&	3.2575E-06&	5.179	&&1.711E-04	&5.257E-04&	1.40\\
3		&1.6269E-06&	2.8366E-06&	7.800	&1.6267E-06&	2.8362E-06&	7.769	&&1.741E-05	&4.346E-05&	2.50\\
5		&3.1532E-06&	4.2934E-06&	13.182	&3.1531E-06&	4.2933E-06&	13.104	&&1.842E-04	&5.694E-04&	4.10\\
7		&3.5801E-06&	5.1314E-06&	18.19	&3.5799E-06&	5.1312E-06&	18.018	&&1.376E-04	&4.759E-04&	6.00\\
10		&3.3621E-06&	4.8638E-06&	20.748	&3.3620E-06&	4.8635E-06&	26.198	&&1.1691E-04&4.1396E-04&	8.80\\
\bottomrule
\end{tabular}
\end{table}

\begin{table}[!htbp]
\caption{Comparison of the mExp-DQM solutions of Example \ref{ex5} with $h=0.1,\delta=0.01,\alpha=1,\beta=1,p=1$}\label{tab-5.1}
\vspace{.2cm}
\centering
\begin{tabular}{l*7l}
\toprule
$t$&  \multicolumn{3}{l}{ mExp-DQM}  &{}& \multicolumn{3}{l}{ Mittal and Bhatia \cite{MB14}}   \\ \cline{2-4} \cline{6-8}
{} &$L_2$ &$L_\infty$ & CPU(s)&{}&$L_2$ &$L_\infty$ &CPU(s)\\
\midrule
1	$\qquad\quad$&3.9796E-04$\qquad\quad$&	6.7076E-04$\qquad$&	0.031$\qquad$&&	1.4441E-02$\quad\qquad$&	2.9996E-02$\qquad$&	0.03\\
2	&4.5099E-05&	1.1091E-04&	0.063&&	1.3898E-03&	3.9711E-03&	0.05\\
3	&4.0589E-05&	7.4545E-05&	0.109&&	1.3018E-03&	2.2178E-03&	0.08\\
5	&4.1078E-06&	8.6460E-06&	0.187&&	1.1112E-04&	2.0618E-04&	0.11\\
7	&4.6749E-07&	1.0452E-06&	0.234&&	1.3695E-05&	3.0052E-05&	0.14\\
10	&3.8692E-08&	7.0454E-08&	0.312&&	1.4408E-06&	2.5354E-06&	0.19\\
\bottomrule
\end{tabular}
\end{table}

\begin{table}[!htbp]
\caption{Comparison of the mExp-DQM solutions of Example \ref{ex5} with $h=.05,\delta=0.001,\alpha=1,\beta=1,p=1$}\label{tab-5.2}
\vspace{.2cm}
\centering
\begin{tabular}{l*8lr}
\toprule
$t$&  \multicolumn{3}{l}{ mExp-DQM}  &{}& \multicolumn{3}{l}{ Mittal and Bhatia \cite{MB14}}   \\ \cline{2-4} \cline{6-8}
{} &$L_2$ &$L_\infty$ & CPU(s)&{}&$L_2$ &$L_\infty$ &CPU(s)\\
\midrule
0.5	$\qquad\quad$&1.28E-04$\qquad\quad$&	2.67E-04$\qquad$&	1.045$\qquad$&&		3.4808E-03$\quad\qquad$&	9.5129E-03$\qquad$&	0.50\\
1	&1.05E-04&	1.82E-04&	2.074 &&	3.2351E-03&	7.4749E-03&	0.70\\
2	&1.04E-05&	3.07E-05&	4.181 &&	2.8518E-04&	1.0361E-03&	1.30\\
3	&1.09E-05&	2.05E-05&	6.255 &&	3.1028E-04&	5.7859E-04&	1.90\\
5	&1.03E-06&	2.40E-06&	10.358&&	2.4495E-05&	6.7234E-05&	3.30\\
7	&1.03E-07&	2.59E-07&	14.446&&	2.5376E-06&	8.2203E-06&	3.90\\
10	&1.18E-08&	2.18E-08&	20.498&&	3.6505E-06&	8.5897E-06&	5.20\\
\bottomrule
\end{tabular}
\end{table}

\begin{table}[!htbp]
\caption{Comparison of the mExp-DQM solutions of Example \ref{ex6} with $h=0.1,\delta=0.01,\alpha=\beta=1$}\label{tab-ex6.1}
\vspace{.2cm}
\centering
\begin{tabular}{lll*7l}
\toprule
$t$ & \multicolumn{5}{l}{mExp-DQM ($p=0.5, 1$)}  &$\quad$ & \multicolumn{3}{l}{ Mittal and Bhatia \cite{MB14} }   \\ \cline{2-6} \cline{8-10}
{}	& $L_2 : p=1$	& $L_\infty: p=1$	&	    $L_2: p=0.5$ &	$L_\infty: p=0.5$ &	CPU(s)	& & $L_2$	& $L_\infty$	& CPU(s) \\
\midrule										

1	$\quad$&5.7365E-04	$\qquad$&7.1586E-04		$\quad$&5.7365E-04	$\quad$&7.1591E-04	$\quad$&0.05	&&1.6144E-03$\qquad$&	3.6006E-03$\quad$&	0.07\\
2	&1.7371E-04	&2.2392E-04		&1.7372E-04	&2.2396E-04	&0.11	&&2.6345E-03&	5.7068E-03&	0.09\\
3	&1.9296E-05	&2.1468E-05		&1.9298E-05	&2.1476E-05	&0.14	&&5.3845E-04&	1.2479E-03&	0.11\\
5	&6.4893E-06	&8.5658E-06		&6.4900E-06	&8.5675E-06	&0.25	&&1.2418E-04&	2.1003E-04&	0.15\\
7	&1.3028E-06	&1.6270E-06		&1.3028E-06	&1.6272E-06	&0.33	&&1.3653E-05&	2.6261E-05&	0.12\\
10	&5.8266E-08	&7.3567E-08		&5.8266E-08	&7.3573E-08	&0.45	&&7.5592E-06&	1.4083E-06&	0.20\\
\bottomrule
\end{tabular}
\end{table}

\begin{table}[!htbp]
\caption{Comparison of the mExp-DQM solutions of Example \ref{ex6} with $h=0.05,\delta=0.001,\alpha=\beta=1, p=0.5, 1$}\label{tab-ex 6.2}
\vspace{.2cm}
\centering
\begin{tabular}{llll*6l}
\toprule
$t$&  \multicolumn{5}{l}{   mExp-DDQM (p=0.5, 1)  }  &\qquad& \multicolumn{3}{l}{ Mittal and Bhatia \cite{MB14}}  \\ \cline{2-6} \cline{8-10}
{}	& $L_2 : p=1$	& $L_\infty: p=1$	&	    $L_2: p=0.5$ &	$L_\infty: p=0.5$ &	CPU(s)	& & $L_2$	& $L_\infty$	& CPU(s) \\
\midrule
0.5	$\quad$&3.2617E-05$\quad$&	4.6301E-05$\quad$&		3.2617E-05$\quad$&	4.6306E-05$\quad$&	1.33$\quad\qquad$&&	3.5833E-04$\quad$&	9.5129E-04$\quad$&	 0.3\\
1	&5.5100E-05&	7.2237E-05&		5.5100E-05&	7.2239E-05&	2.54&&	3.2351E-04&	7.4749E-04&	0.7\\
2	&1.5539E-05&	2.1379E-05&		1.5539E-05&	2.1380E-05&	5.09&&	2.8518E-05&	1.0361E-04&	1.3\\
3	&8.3598E-07&	1.1501E-06&		8.3602E-07&	1.1504E-06&	7.61&&	3.1028E-05&	5.7859E-04&	1.7\\
5	&5.1811E-07&	7.4769E-07&	    5.1812E-07&	7.4772E-07&	12.81&&	2.4495E-06&	6.7234E-05&	2.9\\
7	&1.5582E-07&	1.9983E-07&	 	1.5582E-07&	1.9983E-07&	17.78&&	2.5376E-07&	8.2203E-07&	4.1\\
10	&7.1281E-09&	9.0990E-09&  	7.1281E-09&	9.0991E-09&	25.80&&	3.6505E-09&	8.5897E-08&	5.4\\
\bottomrule
\end{tabular}
\end{table}

\begin{table}[!htbp]
\caption{Comparison of the mExp-DQM solutions of Example \ref{ex7} with $\alpha=\beta=1,\Delta t=0.001, p= 1$ and $h=0.05$}\label{tab-ex7.0}
\vspace{.2cm}
\centering
\begin{tabular}{*9c*6c}
\toprule
$t$&  \multicolumn{3}{l}{mExp-DQM }  &{}& \multicolumn{4}{l}{Mittal and Bhatia \cite{MB14}} &{}& \multicolumn{4}{l}{ Dehghan and Ghesmati \cite{DG10a}}  \\ \cline{2-5} \cline{7-10} \cline{12-15}
{} &	$L_2$ & $L_\infty$ &	$R_e$ &	CPU (s) &&	$L_2$ &	$L_\infty$ &	$R_e$ &	 CPU (s) &&	$R_e:$MLWS &	CPU (s) &	$R_e:$MLPG	& CPU (s) \\	
\midrule
0.5	&4.795E-05	~~~&9.727E-05	~~&1.097E-03	&1.05	&&1.069E-03	~~~&2.474E-03	~~&1.109E-03	&0.5	&&7.939E-05	&9.2	&9.991E-05	&21.0\\
1	&7.290E-05	&1.081E-04	&1.394E-03	&2.18	&&1.529E-03	&3.308E-03	&1.327E-03	&1.1	&&9.098E-05	&12.9	&7.198E-05	&36.2\\
2	&2.946E-05	&4.931E-05	&4.466E-04	&4.23	&&4.647E-04	&1.138E-03	&3.195E-04	&2.0	&&8.705E-04	&25.7	&8.784E-05	&49.1\\
3	&1.200E-05	&2.086E-05	&1.567E-04	&6.28	&&2.199E-04	&4.358E-04	&1.302E-04	&2.8	&&9.931E-04	&38.1	&4.801E-04	&66.8\\
4	&1.281E-05	&1.948E-05	&1.502E-04	&8.46	&&2.715E-04	&5.414E-04	&1.444E-05	&4.3	&&4.703E-03	&49.8	&6.091E-04	&82.0\\
5	&7.989E-06	&1.247E-05	&8.626E-05	&10.49	&&1.720E-04	&3.481E-04	&8.423E-05	&7.0	&&7.302E-03	&62.0   &9.498E-04	&97.3\\
10  &2.738E-06  &	4.198E-06&	2.314E-05&	21.00&&	7.729E-05&	1.404E-04&	 2.962E-05&	9.6	&&		  &	    &            &      \\
\bottomrule
\end{tabular}
\end{table}



\end{document}